\documentclass[11pt,twoside]{article}
\usepackage{amssymb}
\usepackage[mathscr]{eucal}
\usepackage{amsmath}
\usepackage{mathrsfs}
\usepackage{braket}
\usepackage{xcolor}

\def\chiu{\hfill$\displaystyle\vspace{4pt}
\underset\Box\null$}
\def\Pr{{\bf Proof. }}
\def\O{\Omega}

\def\R{\Bbb R}
\def\N{\Bbb N}
\def\o{\"{o}}
\def\à{\`{a}}
\def\è{\`{e}}
\def\ì{\`{i}}
\def\ù{\`{u}}
\def\ò{\`{o}}
\def\é{\'{e}}
\def\vf{\varphi}
\def\dy{\displaystyle}

\def\be{\begin{equation}}
\def\ba{\begin{array}}
\def\ea{\end{array}}
\def\ee{\end{equation}}
\def\vs1{\vspace{1ex}}
\def\cD{{\mathcal D}}
\def\vp{\varphi}

\def\po{\partial\Omega}

\font\sc=cmcsc10
\title{Existence of global weak solutions to a parabolic $p$-Laplacian problem with convective term}
\author{\sc Angelica Pia Di Feola\thanks{
Dipartimento di Matematica e Fisica, 
Universit\`{a} degli Studi della Campania ``Luigi Vanvitelli'', viale Lincoln 5, 81100 Caserta,
 Italy. E-mail: angelicapia.difeola@unicampania.it}
$\;$  and Michael R{\r u}{\v z}i{\v c}ka\thanks{Institute of Applied Mathematics, Albert-Ludwigs-University, Ernst-Zermelo-Straße 1, 79104 Freiburg, Germany.
E-mail: rose@mathematik.uni-freiburg.de}}
 \date{}
\begin{document}
\maketitle
\noindent{\bf Abstract}  - {\small For a given
    bounded domain $\O \subset \R^3$, with $C^2$  boundary, and a
    given instant of time $T>0$, we prove the existence of a global
  weak solution on $(0,T)$, which satisfies a maximum
  principle, to a parabolic $p$-Laplacian system with convective term
  without divergence constraint for any $p\in (1,2)$. }\vskip 0.2cm
 \vskip -0.7true cm\noindent
\newcommand{\red}{\protect\bf}
\renewcommand\refname{\centerline
{\red {\normalsize \bf References}}}
\numberwithin{equation}{section}
\newtheorem{tho}{\bf Theorem}[section]

\newtheorem{lemma}[tho]{\bf Lemma}
\newtheorem{prop}[tho]{\bf Proposition}
\newtheorem{coro}[tho]{\bf Corollary}
\newtheorem{ass}[tho]{\bf Assumption}
\newtheorem{defi}[tho]{\bf Definition}
\newtheorem{rem}[tho]{\sc Remark}
\renewcommand{\theequation}{\thesection .
\arabic{equation}}
\setcounter{section}{1}
\section*{Introduction}
\renewcommand{\theequation}{1.\arabic{equation}}
\renewcommand{\thetho}{1.\arabic{tho}}
\renewcommand{\thedefi}{1.\arabic{defi}}
\renewcommand{\therem}{1.\arabic{rem}}
\renewcommand{\theprop}{1.\arabic{prop}}
\renewcommand{\thelemma}{1.\arabic{lemma}}
\setcounter{tho}{0} 
This note deals with the $p$-Laplacian system with convective term,
in the subquadratic case $p<2$:
\begin{align}
\partial_t v-\nabla\cdot\left(|\nabla v|^{p-2}\nabla v\right) &=  v\cdot\nabla v\,,& &\textrm{
in }(0,T)\times\O, \nonumber 
\\
v(t,x)\dy &=0\,, &   &\textrm{ on }
(0,T)\times\po, \label{PF}
\\
v(0,x)&=v_\circ(x),&&\mbox{ on
}\{0\}\times\O, \nonumber
\end{align} where  $\O$ is a given bounded
domain in $\R^3$, whose boundary is $C^2$-smooth, ${T>0}$ is a given
instant of time, $v_\circ$ is the initial datum belonging to $
L^{\infty}(\O)$, and ${v:(0,T)\times\O\to \R^3}
$ is the sought vector field.
\par The interest in the system \eqref{PF}, already studied in
\cite{CDF}, comes from its relation to the $p$-Navier-Stokes
equations, which model the behavior of so-called power-law
fluids. The motion of such fluids is described by
\begin{align}
 \partial_t u-\nabla\cdot \left(|\cD u |^{p-2}\cD u\right)+
\nabla\pi & = - u\cdot\nabla u+f,& \textrm{ in
}(0,T)\times\O,
\label{PLF} \\
\nabla\cdot u&=0, & \textrm{ in }(0,T)\times \O, \nonumber
\end{align}  
where $u$ denotes the velocity field, $\pi$ the pressure, $f$ is a given external force, and $\cD =\frac 12( \nabla u +\nabla u ^\top)$ is the symmeytric gradient. For background on both the physical modeling and the mathematical theory of non-Newtonian fluids, we refer the reader to \cite{15, MNRR,  Raj, show1}.
System \eqref{PF} is obtained from \eqref{PLF} by
dropping the divergence constraint and the related pressure gradient, and by replacing the symmetric gradient with the full gradient. Clearly, these  changes make the system  not suitable as a model for fluid dynamics, as the corresponding {\it constitutive law } becomes not compatible with the principle of material
invariance. 
\par  It is well known that for the Navier-Stokes equations, obtained by
setting $p = 2$  in \eqref{PLF}, if one examines the linearized Stokes
problem (by dropping the convective term) one has the existence of
global regular solutions (see  \cite{Sol}). Furthermore, for the
unsteady Stokes problem the validity of a maximum modulus theorem was proved in \cite{AG} for bounded domains and in \cite{Abe, M, M2} for exterior domains. 
 On the other hand, if one considers the purely parabolic problem associated to the Navier-Stokes equations, obtained by eliminating the pressure term and the solenoidal constraint, both existence of global regular solutions and a maximum principle hold. For details on these results, we refer to \cite{Avrin,Mar}.
 \par The purpose of this paper is to continue the investigations of the purely parabolic problem \eqref{PF}, that is the counterpart of the above problems when the power-law fluids system \eqref{PLF} is considered in place of Navier-Stokes one.
The first work in this direction is \cite{CDF}:
   apart from it, we are not aware of any other existence results or
   $L^\infty$-estimates for system \eqref{PF}. Here the existence of
   strong solutions and the maximum principle is proven for $p\in
   \left(\frac{3}{2}, 2\right)$ and an initial datum ${u_\circ \in
   L^\infty(\O)}
   $. The assumption $p>\frac{3}{2}$ is related to the
   regularity of the solution. In fact, it is required in order to
   obtain the $ L^{2}(0,T; L^\frac{4}{4-p}(\O))$ estimate of weighted
   second derivatives  by means of Lemma 2.4 in \cite{CDF}. 
In the present paper,  we investigate the existence of weak solutions  for $(\ref{PF})$ for $p\in (1,2)$ and the innovation  lies precisely in extending the range of $p$. Indeed, we succeed in establishing the existence of global weak solutions for every $p\in (1,2)$.
We construct a weak solution of $(\ref{PF})$ in sense of Definition
$\ref{defnomu}$ through a procedure that does not rely on uniform
second order derivative estimates and the consequent strong convergence of gradients. Our approach is as follows: we study an approximating system similar to the one in \cite{CDF}, establish uniform estimates in $L^\infty(0,T;L^2(\O)) \cap L^p(0,T;W^{1,p}_0(\O))$ and the maximum principle for the sequence of approximations, and apply Minty’s trick in the passage to the limit.
Concerning the maximum principle, as in \cite{CDF}, the idea highly relies on that in \cite{Mar},  where, for the first time to the best of our knowledge, the duality technique has been employed for proving the maximum principle for a quasilinear  system.
A crucial novelty for our approach is an integration by parts formula (cf.~Proposition $\ref{ibp}$) for the following space 
$$W 
\Set{\!  \! u \! \in \! L^{\infty}(\Omega_T) \! \cap\! L^{p}(0,T; W^{1,p}_0(\Omega)) \!|\! \partial_t u \!\in \! L^{p}(0,T; L^{p}(\Omega)) \! + \! L^{p'}(0,T;W^{-1,p'}(\Omega))\!\!} \!,$$
which, to the best of our knowledge, is not known before. 
The approach is the same as in \cite{D}: at first, we establish a
suitable density result (cf.~Proposition~\ref{den}) and then prove a
continuity result (cf.~Proposition \ref{con}).
\par
To formulate our main result we define the notion of solutions we are
working with: 
  \begin{defi}\label{defnomu}
{\rm Let  $v_{\circ}\in L^2(\O)$. A vector field
$v\!:(0,T)\times \O\to\R^3$
 is said to be a global weak solution of the 
system {\rm \eqref{PF}} if  
$ v\in 
L^{p}(0,T; W_0^{1,p}(\O))\cap L^{p'}(0,T;L^{p'}(\O))$ satisfies 
$$\int_0^T\!\!\! \left[(v, \partial_\tau \psi)-\left(|\nabla v|^{p-2}\,
\nabla v,\nabla \psi\right)-(v\cdot \nabla v, \psi)\right] d\tau=-(v_\circ, \psi(0)),   $$
for all $ \psi\in C_0^\infty([0,T)\times \O)$, and
$$\lim_{t\to 0^+}(v(t), \vp)=(v_\circ, \vp),$$
for all $ \vp \in C_0^\infty(\O)\,.$}
\end{defi}
Our main result is the following existence theorem of global weak solutions corresponding to an initial datum in $L^{\infty}(\O)$, together with a maximum principle: 
\begin{tho}\label{exis} 
{\sl Assume that $v_\circ$
belongs to $L^{\infty}(\O)$. For all $p\in (1,2)$ and $T>0$, there exists a global weak solution $v$
of the system $(\ref{PF})$ in the sense of Definition\;$\ref{defnomu}$.
In particular, 
$$\dy \|v\|^2_{{C([0,T];L^2(\O))}}+\|\nabla v\|^p_{L^{p}(0,T;L^p(\O))}\leq K_1(\|v_\circ\|_\infty,T)\,.$$
Moreover, $v$ satisfies the maximum principle
\begin{equation}\label{mpf}
||v(t)||_{L^{\infty}(\Omega)} \leq ||v_\circ||_{L^{\infty} (\Omega)}, \;\; \; \text{for all }\, t \in [0,T].
\end{equation}
Finally,  $v\in C([0,T]; L^q(\O))$, for any $q\in [1,\infty)$, and
\begin{equation}\label{id}
\lim_{t \to 0^+}||v(t) - v_\circ||_{L^q(\Omega)} =0, \;\; \; \text{for all }\, q \in [1,+\infty).
\end{equation}}
\end{tho}
\par {\it Outline of the proof} -  We introduce a  non-singular
approximating problem with a Laplacian term and a modified convective
term (see \eqref{PFepv}). Through the Bochner  pseudo-monotone
operator theory, we obtain a local existence result and, via a duality
technique, establish the maximum principle for the approximating
solution. Subsequently, by the maximum principle, we get the global
existence of solutions. Finally, using estimates independent of the approximating parameters, we deduce the validity of the existence result and the maximum principle for a solution of \eqref{PF}, constructed as the limit of the approximations.
\section{\normalsize  Notations and auxiliary results}\renewcommand{\theequation}{2.\arabic{equation}}
\renewcommand{\thetho}{2.\arabic{tho}}
\renewcommand{\thedefi}{2.\arabic{defi}}
\renewcommand{\therem}{2.\arabic{rem}}
\renewcommand{\theprop}{2.\arabic{prop}}
\renewcommand{\thelemma}{2.\arabic{lemma}}
\renewcommand{\thecoro}{2.\arabic{coro}}
\setcounter{tho}{0}
Throughout the paper, the constants may depend on $T$ and $\O$. Moreover, their value may change from line to line.
We work within the framework of standard Lebesgue spaces
$L^p(\Omega)$, $p\in [1,\infty]$, and Sobolev spaces
$W^{k,p}(\Omega)$, ${p\in[1,\infty], k\in\mathbb
  N}$. We use the notation $(f,g):= \dy \int_\O fg \,
  dx $, whenever it is well defined, {and $\langle \cdot, \cdot \rangle_X$ for the duality pairing between a Banach space $X$ and its dual $X'$}. In addition, we consider the classical Bochner spaces $L^p(0,T;X)$, where $X$ is a Banach space, $p \in [1, \infty]$, and frequently make use of the identity $$L^q(0,T; L^q(\O))=L^q(\O_T), \text{ for } q\in[1,\infty),$$ where  $\O_T:=(0,T)\times \O$. For more details we refer to \cite{ABS}.
We recall a well-known result that will be used in our analysis.
\begin{tho}\label{aubin}{\rm (Aubin-Lions)} {\sl Let $X$, $X_1$, $X_2$ be Banach spaces. Assume that   $X_1$ is compactly embedded in $X$ and $X$ is continuously embedded in $X_2$, and that $X_1$ and $X_2$ are reflexive.
 For $1<q,s<\infty$, set $$W=\{\psi \in L^s(0,T; X_1)\,\big |\,
 \partial _t\psi\in L^q(0,T; X_2)\}\,.$$
Then the inclusion $W \subset L^s(0,T; X)$ is compact.} \end{tho}
In the sequel we use some properties of mollifiers (cf.~\cite{ABS, RRS} for more details).
\begin{lemma}
Let $(\rho_\varepsilon)_{\varepsilon > 0} \subset C^\infty(\R^n)$ be a standard family of spatial mollification kernels.
Then, for any function $v$ in $L^{p}_{\text{loc}}(\mathbb{R}^n)$ with $1 \leq p < \infty$, the convolution
\[
J_\varepsilon(v)(x) := (\rho_\varepsilon * v)(x) = \int_{\mathbb{R}^n} \rho_\varepsilon(y)\, v(x - y)\, dy
\]
defines a smooth approximation of $v$ with the following properties:
\begin{itemize}
 \item $J_\varepsilon v \in C^\infty(\mathbb{R}^n)$, for all $\varepsilon > 0$;
 \item if $v \in L^{p}(\O)$ where $1\leq p < \infty$, then $J_\varepsilon (v) \in L^{p}(\O)$ with
\begin{equation}
\label{M1} ||J_\varepsilon ( v )||_{L^{p}(\Omega)} \leq ||v||_{L^{p}(\Omega)}  \; \; \; \text{ and } \; \; \; \lim_{\varepsilon \to 0^+} ||J_\varepsilon (v) - v||_{L^{p}(\Omega)}  =0;
\end{equation}
 \item if $v \in L^{\infty}(\O)$, then $J_\varepsilon ( v )\in L^{\infty} (\O)$ with
\begin{equation}
\label{M3} ||J_\varepsilon ( v )||_{L^{\infty}(\O)} \leq ||v||_{L^{\infty}(\O)};
\end{equation}
\item if $v \in L^2(\Omega )$, then
\begin{equation}
\label{M2} ||J_\varepsilon (v) ||_{L^{\infty}(\O)} \leq c(\varepsilon) ||v||_{L^2(\Omega )}
\end{equation}
with $c(\varepsilon) = \varepsilon^{-\frac{n}{2}}\to + \infty$ as $\varepsilon \to 0^+.$
\end{itemize}
\end{lemma}
 \begin{lemma} Let $f \in L^p(\R;X)$ for some $p \in [1,\infty)$, and $\varphi \in L^1(\R)$. For $\varepsilon > 0$, denote $\varphi_\varepsilon(y) := \varepsilon^{-d} \varphi(\varepsilon^{-1}y).$ Then $$\varphi_\varepsilon * f \;\to\; c_\varphi f \quad \text{in } L^p(\R;X)$$
as $\varepsilon \to 0^+$, where $c_\varphi := \int_{\R} \varphi(y)\, dy.$ \end{lemma}
For  a suitable given $v$,  we denote by $\widetilde
J_{\eta}(\nabla v)$ the convolution $\widetilde \rho_\eta * \nabla v$,
where $\widetilde \rho _\eta$ is a time-spatial  mollifier, i.e., $\rho $ is a standard mollification kernel in
$\R^{3+1}$, and for  $p\in (1,2)$, $\mu, \eta>0$, we set for $(t,x)\in (0,T)\times \O$
\be \label{amu} a_\eta(\mu,v)(t,x):=\left(\mu+|\widetilde J_{\eta}(\nabla
v)(t,x)|^2\right)^\frac{p-2}{2},\ee
and
$$ a(\mu,v)(t,x):=\left(\mu+|\nabla
v(t,x)|^2\right)^\frac{p-2}{2}.$$
 \begin{lemma}\label{comp}{\sl Assume $\mu>0$, and $t>0.$
 Assume that $\nabla\psi\in L^2(0,t; L^{2}(\Omega))$, $\nabla v \in L^r(0,t;L^r(
\O))$, for some $r>1$, and let $h^m$ be a
 sequence such that $\nabla h^{m}$ is uniformly bounded in $ L^2(0,t;
 L^{2}(\Omega))$ with respect to $m\in \N$. Then, there exists a subsequence $h^{m_k}$ such
 that
 $$\lim_{k\to \infty}\int_0^t \int_\O\left(a_{\frac{1}{m_k}}(\mu,v)- a(\mu,v)\right)
\nabla h^{m_k}\cdot  \nabla
\psi \,dx\,d\tau=0\,.$$ }
 \end{lemma}
\Pr
The proof works as in \cite[Lemma 2.5]{CMparabolic}. For completeness, we reproduce it.
Since $\nabla v \in L^r(0,t;L^r(\O))$, $\widetilde J_{\frac {1}{m}}(\nabla v) \to \nabla v$ in $L^r(0,t;L^r(\O))$, so there exists a subsequence such that  $\widetilde J_{\frac {1}{m_k}}(\nabla v) (s,x) \to \nabla v (s,x)$ almost everywhere in $(0,t) \times \O$.
\\By Hölder's inequality
\begin{align*}
&\int_0^t \int_\O\left(a_{\frac{1}{m_k}}(\mu,v)- a(\mu,v)\right)
\nabla h^{m_k}\cdot  \nabla
\psi \,dx\,d\tau\\ 
&\leq \left( \int_0^t \int_\O\left|a_{\frac{1}{m_k}}(\mu,v)- a(\mu,v)\right|^2
\left| \nabla \psi \right|^2 \,dx\,d\tau\right)^\frac{1}{2} || \nabla h^{m_k} ||_2.\end{align*}
Moreover,
$$ a_{\frac{1}{m_k}}(\mu,v)(x,t)-
 a(\mu,v)(x,t) \to 0 \text{ a.e. in }(0,t)\times\O,$$
 and 
 $$\left|a_{\frac{1}{m_k}}(\mu,v) - a(\mu,v) \right|^2
\left| \nabla \psi \right|^2 \leq \left( 2\mu^{\frac{p-2}{2}}
\right)^2 \left| \nabla \psi \right|^2. %
$$%
By Lebesgue dominated convergence, since $\left| \nabla \psi \right|^2 \in L^1(0,t; L^1(\O)),$ we get the assertion.
\chiu

We recall some results from the theory of Bochner pseudo-monotone
operators, needed for our purposes (cf.~\cite{KR}).
\begin{defi}[\rm Bochner pseudo-monotonicity]\label{psemon} Let $T \in
  (0,\infty)$, and $1 < p < \infty$.
Let $(V,H,j)$ be a Gelfand-Tripel, where $V$ is a reflexive, separable Banach space. 
An operator $\mathcal{A}: D(\mathcal{A}) \subset L^p(0,T;V)\cap_j L^\infty(0,T;H) \to L^{p'}(0,T;V')$ 
is said to be Bochner pseudo-monotone if from
\begin{align*}
    u_n \rightharpoonup u \quad  &\text{in } L^p(0,T;V) \quad (n \to \infty), \\
    j u_n \overset{\ast}{\rightharpoonup} ju \quad &\text{in } L^\infty(0,T;H) \quad (n \to \infty),\\
    j u_n(t) \rightharpoonup j u(t) \quad  &\text{in } H \text{ for a.e. } t\in 0,T \quad (n \to \infty),
    \end{align*}
and
$$ \limsup_{n \to \infty} \langle\mathcal{A}u_n, u_n - u\rangle _{L^p(0,T;V)}  \leq 0, $$
it follows that
$$ \langle \mathcal{A}u, u-w \rangle _{L^p(0,T;V)} 
\leq \liminf_{n \to \infty} \langle \mathcal{A}u_n, u_n - w \rangle _{L^p(0,T;V)},$$
for all  $w \in L^p(0,T;V).$
\end{defi}
\begin{defi}[\rm Bochner coercivity]\label{bocoer} Let $T \in
  (0,\infty)$, and $1 < p < \infty$. 
Let $(V,H,j)$ be a Gelfand-Tripel, with a reflexive, separable Banach space $V$.
An operator $\mathcal{A}: D(\mathcal{A}) \subset L^p(0,T;V)\cap_j L^\infty(0,T;H) \to L^{p'}(0,T;V')$ 
is said to be Bochner coercive w.r.t. $f \in L^{p'}(0,T;V')$ and $v_\circ \in H$ if  there exists a constant $M:=M(f, v_\circ, \mathcal{A})> 0$ such that for all $v \in L^p(0,T;V)\cap_j L^\infty(0,T; H)$ from
$$\frac{1}{2} ||jv(t)||_H^2 + \langle\mathcal{A}v-f, x\chi_{[0,T]}\rangle _{L^p(0,T;V)} \leq \frac{1}{2} ||v_\circ||_H^2 \quad \text{ for a.e. } t \in (0,T),$$
it follows that $||x||_{L^p(0,T;V)\cap_j L^\infty(0,T;H)} \leq M$. 
Moreover, it is said to be Bochner coercive if it is Bochner coercive with respect to $f$ and $v_\circ$ for all $f\in L^{p'}(0,T;V')$ and $ v_\circ \in H$.
\end{defi}
For the reader’s convenience, we recall some definitions and results established in \cite{KR} in the setting relevant for our purposes:
$(V,H,j) = (L^2(\Omega), W^{1,2}_0(\Omega),  \mathrm{id})$ and $p=2$. 
For the space 
$$\mathcal{W}:= \set{ v \in L^2(0,T; W^{1,2}_0(\O)) \,\big |\, \partial_t v \in L^2(0,T; W^{-1,2}(\O))}$$
the following result is proved in \cite[Proposition 2.23]{KR}.
\begin{prop}\label{contrap}
For $0<T<+\infty$, it holds:
\begin{itemize}
    \item[(i)] the space $\mathcal{W}$ forms a Banach space equipped with the norm
    $$||v ||_{\mathcal{W}}:= ||v ||_{L^2(0,T; W^{1,2}_0(\O))} + ||\partial_t v  ||_{L^2(0,T; W^{-1,2}(\O))}.$$
    \item[(ii)] Given $v \in \mathcal{W}$, it possesses a unique representation in $C([0,T]; L^2(\O))$  and the resulting mapping is an embedding.
    \item[(iii)] The following generalized integration by parts formula holds:
\begin{align} \label{intbyp}   
\int_{t'}^t \! \langle \partial_s u(s), v(s)  \rangle _{ W^{1,2}_0(\O)} ds &=  (u(t), v(t)) - (u(t'), v(t'))\nonumber\\
&\quad - \int_{t'}^t \! \! \langle  u(s), \partial_s v(s) \rangle _{ W^{1,2}_0(\O)} ds 
\end{align}
for all $u,v \in \mathcal{W}$, and $t,t' \in [0,T]$ with $t'\leq t.$
\end{itemize}
\end{prop}
\begin{defi}[\rm Induced operator]
   Let  $0<T<+\infty$. We call an operator $\mathcal{A}: L^2(0,T;W^{1,2}_0(\O))\cap L^\infty(0,T;L^2(\O)) \to L^2(0,T;W^{-1,2}(\O))$ induced by a family of operators $A(t): W^{1,2}_0(\O)\to W^{-1,2}(\O)$, $t \in (0,T)$, if
    $$\langle\mathcal{A}u, v\rangle _{L^2(0,T;W^{1,2}_0(\O))}:= \int_0^T \langle A(t)(u(t)), v(t)\rangle _{W^{1,2}_0(\O)} \, dt,$$
    for every $u \in L^2(0,T;W^{1,2}_0(\O))\cap L^\infty(0,T;L^2(\O))$, and $v \in L^2(0,T;W^{1,2}_0(\O))$.
\end{defi}
We recall Proposition 3.13 and Theorem 4.1 in \cite{KR}, in the
special case 
$(V,H,j)=(W^{1,2}_0(\O),L^2(\O),id)$, and $p=2$.
\begin{prop}\label{prop13}
Let $A(t): W^{1,2}_0(\O)\to W^{-1,2}(\O)$, $t\in (0,T) \subset \R^+$, be a family of operators such that the following properties (C.1)-(C.5) are fulfilled: 
\begin{itemize}
     \item[(C.1)]  The operator $A(t): W^{1,2}_0(\O) \to W^{.1,2}_0(\O)$ is demi-continuous for almost every $t\in (0,T)$.
     \item[(C.2)] The function $A(\cdot)v : (0,T) \to W^{-1,2}(\O)$ is Bochner-measurable for all ${v\in W^{1,2}_0(\O)}$.
    \item[(C.3)] For some non-negative functions $\alpha, \gamma \in L^2(0,T), \beta \in L^{\infty} (0,T)$ and a non-decreasing function $\mathcal{B}: \R_{\geq 0} \to R_{\geq 0}$ holds
      $$||A(t)v||_{W^{-1,2}(\O)}\leq \mathcal{B} (||v||_{L^2(\O)}) (\alpha(t) + \beta(t) ||v||_{W^{1,2}_0(\O)}) + \gamma (t)$$
     for almost every $t\in (0,T)$, and all $v\in W^{1,2}_0(\O)$.
    \item[(C.4)]  The operator $A(t): W^{1,2}_0(\O) \to W^{-1,2}(\O)$ is pseudo-monotone for almost every $t \in (0,T)$.
    \item[(C.5)] For some constant $c_0> 0$, non-negative functions $c_1, c_2 \in L^1(0,T)$, and a non-decreasing function $\mathcal{C}: \R_{\geq 0} \to \R_{\geq 0} $ holds
    $$\langle A(t)v,v\rangle _{W^{1,2}_0(\O)} \geq c_0 ||v||_{W^{1,2}_0(\O)}^2 - c_1(t)\mathcal{C} (||v||_{L^2(\O)}) - c_2(t)$$
    for almost every $t\in (0,T)$, and all $v\in W^{1,2}_0(\O) $.
\end{itemize}
Then the induced operator $\mathcal{A}: L^2(0,T;W^{1,2}_0(\O))\cap L^\infty(0,T;L^2(\O)) \to L^2(0,T;W^{-1,2}(\O))$ is Bochner pseudo-monotone.
\end{prop}
\begin{tho}\label{theo1}
    Let $A(t): W^{1,2}_0(\O)\to W^{-1,2}(\O)$, $t\in (0,T)
    \subset \R^+$, be a family of operators such that (C.1)-(C.3) are
    fulfilled, and that the induced operator ${\mathcal{A}:
    L^2(0,T;W^{1,2}_0(\O))\cap L^\infty(0,T;L^2(\O)) \to
    L^2(0,T;W^{-1,2}(\O))}$ is Bochner pseudo-mono\-tone and Bochner
    coercive with respect to $f\in L^2(0,T; W^{-1,2}(\O))$ and
    $v_\circ \in L^2(\O)$. Then, there exists a function $v\in
    \mathcal{W}$ such that for all $\vp \in L^2(0,T;W^{1,2}_0(\O))\cap
    L^\infty(0,T;L^2(\O)) $ there holds
    \begin{align}
    \int_0^T \! \! \langle  \partial_t v(t), \vp(t)\rangle _{ W^{1,2}_0(\O)
      }\, dt  + \int_0^T \! \! \langle {A}(t)(v(t)) , \vp(t) \rangle_{ W^{1,2}_0(\O) } \, dt &= \int_0^T \!\! \langle f(t), \vp (t)\rangle _{ W^{1,2}_0(\O) }dt \nonumber\\
    v(0)&=v_\circ . \label{evsy}
    \end{align} 
\end{tho}
\section{\normalsize The approximating systems}\renewcommand{\theequation}{3.\arabic{equation}}
\renewcommand{\thetho}{3.\arabic{tho}}
\renewcommand{\thedefi}{3.\arabic{defi}}
\renewcommand{\therem}{3.\arabic{rem}}
\renewcommand{\theprop}{3.\arabic{prop}}
\renewcommand{\thelemma}{3.\arabic{lemma}}
\renewcommand{\thecoro}{3.\arabic{coro}}
\setcounter{tho}{0}
We  obtain a weak solution of (\ref{PF}) as the limit of a sequence of  weak solutions for the following approximating system
 \begin{align} \partial_t v-\nu\Delta
v-\nabla\cdot\left((\mu+|\nabla
v|^2)^\frac{p-2}{2}\nabla v\right) &= -J_\mu(v)\cdot \nabla v, & &\textrm{ in }(0,T)\times\O, \nonumber \\
v(t,x) & = 0, & &\textrm{ on }
(0,T)\times\po, \nonumber\\
v(0,x) & =v_\circ(x), & &\mbox{ on
}\{0\}\times\O , \label{PFepv}
\end{align}
where  $ \mu, \nu > 0$, $J_\mu$ is a spatial mollifier, and  $v_o(x) \in
L^{\infty}(\Omega)$. 
The term $\nu \Delta v$ is crucial for employing the duality technique in order to prove the maximum principle. The parameter $\mu$ in the $p$-Laplacian is introduced to avoid, in a first step, the singularity of the operator. 
The spatial mollifier $J_\mu(v)$ plays a key role in the initial step, the proof of Proposition \ref{existenceL}. Indeed, thanks to property $(\ref{M2})$, it allows us to define the operator $\mathcal{A}v$ (see below) since $J_\mu(v) \in L^{\infty}(\O_T)$ even if $v \in L^{\infty}(0,T; L^{2}(\Omega))$, when the maximum principle is not yet available.
\begin{defi}\label{weakmunualpha}
{\rm Let be $\mu, \nu> 0$. Let $v_{\circ}\in L^{\infty}(\O)$ and $T> 0$. A vector field
$v\!:(0,T)\times\O\to\R^3$
 is said to be a global weak solution of the
system {\rm $(\ref{PFepv})$} if
\begin{itemize}
\item[i)]$ v\in  L^{\infty}(0,T;
 L^{2}(\Omega))\cap L^2(0,T; W^{1,2}_0(\Omega)), \;  \partial_t v\in   L^2(0,T; W^{-1,2}(\Omega))$;
\item[ii)]$v$ satisfies the integral identity  
\begin{align*}
& \int_0^t \left[(v,\partial_\tau \psi)-\nu(\nabla v,\nabla\psi)-\left(a(\mu, v)\,\nabla v,\nabla \psi\right)-(J_\mu(v)\cdot \nabla v, \psi)\right] d\tau\\ 
 &=(v(t),\psi(t))-(v_\circ, \psi(0)), 
\end{align*}
for all  $t\in [0,T]$, for all  $\psi\in C^\infty([0,T]; C_0^\infty(\O));$
   \item[iii)]$\dy\lim_{t\to
0^+}\|v(t)-v_\circ\|_{ L^{2}(\Omega)}=0.$\end{itemize}}
\end{defi}
\begin{rem} \label{rmk1} By a density argument and the regularity of solution $v$ in Definition \ref{weakmunualpha}, we can consider  test  functions $\psi \in L^{\infty}(0,T;
 L^{2}(\Omega))\cap L^2(0,T; W^{1,2}_0(\Omega))$ with $\partial_t  \psi\in   L^2(0,T; W^{-1,2}(\Omega))$.
\end{rem}
\subsection{\normalsize Local existence of weak solution of $(\ref{PFepv})$.}
\begin{prop}\label{existenceL}
{\sl Let be $\nu>0$, and $\mu>0$. Assume that $v_\circ$ belongs to
  $L^{\infty}(\O)$. There exists a weak solution
    $v^{\mu,\nu}$ of  the system $(\ref{PFepv})$ in $(0,T^*)$ for $T^*:= \frac{\nu \mu^{\frac 32}}{2\|v_\circ\|_2^2} $ satisfying  i)--iii) in  Definition\, \ref{weakmunualpha} with $T$ replaced by $T^*$. }
\end{prop}
\Pr
In order to prove the existence of  a local solution of
$(\ref{PFepv})$, we use Theorem~\ref{theo1}. 
For $T^*>0$, we set 
\begin{equation}
    \langle A(v), w\rangle _{W^{1,2}_0(\Omega)} \,  := \nu  (\nabla v,\nabla w) + ((\mu + |\nabla v|^2)^\frac{p-2}{2}\nabla v , \nabla w)+ (J_\mu (v) \cdot \nabla v, w), \label{op}
\end{equation}
and
 $$\langle\mathcal{A}v, w\rangle _{L^2(0,T^*; W^{1,2}_0(\Omega)) \cap
   L^{\infty}(0,T^*;  L^{2}(\Omega))} := \int_0^{T^*} \langle A(v(t)), w(t)\rangle _{W^{1,2}_0(\Omega)} \, dt.$$
 
\textit{Bochner coercivity } - The reason why, at this stage, we have
a local solution only, is that the Bochner-coercivity holds only for
small $T^*$. Indeed, supposing for $v \in L^2(0,T^*; W^{1,2}_0(\Omega))
\cap L^{\infty}(0,T^*; L^{2}(\Omega))$ that 
\begin{align}
& \frac12\|v(t)\|_{ L^{2}(\Omega)}^2+\nu \int_0^t\|\nabla v\|_{ L^{2}(\Omega)}^2 \, d\tau + \int_0^t \|(\mu +|\nabla v|^2)^\frac{p-2}{4}\nabla v\|_{ L^{2}(\Omega)}^2 \, d\tau  \nonumber \\ 
& + \int_0^t ((J_\mu(v)\cdot\nabla v), v) \, d\tau  \leq \frac12\|v_\circ\|_{ L^{2}(\Omega)}^2,   \label{ap2}
\end{align}
then it follows by H\o lder's inequality,
$\|J_\mu(v)\|_{L^{\infty}(\O)}\leq \mu^{-\frac{3}{2}}\|v\|_{
  L^{2}(\Omega)}$ (cf. $(\ref{M2})$), and  $(\ref{ap2})$ that 
\begin{align*}
& \frac12 \|v(t)\|_{ L^{2}(\Omega)}^2 +\nu \int_0^t\|\nabla v\|_{ L^{2}(\Omega)}^2\,  d\tau + \int_0^t\|(\mu+|\nabla
v|^2)^\frac{p-2}{4}\nabla v\|_{ L^{2}(\Omega)}^2 \, d\tau \\
& \leq  \mu^{-\frac 32}\int_0^t\|\nabla v\|_{ L^{2}(\Omega)}\|v\|_{ L^{2}(\Omega)}^2\, d\tau +\frac12\|v_\circ\|_{L^2(\O)}^2.
\end{align*}
By applying Young's inequality on the right-hand side, one gets that
 \begin{align}
& \frac12 \|v(t)\|_{L^2(\O)}^2 +\frac{\nu}{2} \int_0^t\|\nabla v\|_{ L^{2}(\Omega)}^2\, d\tau +\int_0^t\|(\mu+|\nabla v|^2)^\frac{p-2}{4}\nabla v\|_{ L^{2}(\Omega)}^2 \, d\tau  \nonumber \\ 
& \leq \frac{\mu^{-\frac 32}}{2\nu}\int_0^t\|v\|_{ L^{2}(\Omega)}^4\, d\tau +\frac12\|v_\circ\|_{ L^{2}(\Omega)}^2. \label{ap2abb} 
\end{align}
We end up with the following integral inequality
\be\label{ap2abc}
\frac12\|v(t)\|_{ L^{2}(\Omega)}^2 \leq \frac{1}{2 \nu \mu^{\frac{3}{2}}}\int_0^t\|v\|_{ L^{2}(\Omega)}^4\, d\tau + \frac12\|v_\circ\|_{ L^{2}(\Omega)}^2,\ee
that  yields by integration with respect to time 
\be\label{ap2ac}
\|v(t) \|_{L^{2}(\Omega)}^2
\leq \frac{\nu \mu^{\frac{3}{2}}\|v_\circ\|_{ L^{2}(\Omega)}^2}{\nu \mu^{\frac{3}{2}}-\|v_\circ\|_{ L^{2}(\Omega)}^2  t},\ \text{ for all } t\in\left[0, \frac{\nu \mu^{\frac{3}{2}}}{\|v_\circ\|_{ L^{2}(\Omega)}^2}\right).\ee
Thus, for $T^* :=\frac{\nu \mu^{\frac{3}{2}}}{2\|v_\circ\|_{ L^{2}(\Omega)}^2}$,
using \eqref{ap2ac}, estimate $(\ref{ap2abb})$ ensures 
\be\label{nunu}\|
v\|_{L^2(0,T^*; W^{1,2}_0(\Omega)) \cap L^{\infty} (0,T^*;  L^{2}(\Omega))}\leq c\,\|v_\circ\|_{ L^{2}(\Omega)}^2\,.\ee
Ultimately, the operator $\mathcal{A}$, on $(0,T^*)$ with $T^*$ defined as above, is Bochner coercive w.r.t. $f=0$ and $v_\circ.$\\

\textit{Bochner pseudo-monotonicity} - We use Proposition \ref{prop13}
and have to check conditions $(C.1)$-$(C.5)$. 
 \begin{itemize}
 \item[(C.1)] The first term in \eqref{op} is a bounded linear operator, so it is demi-continuous. The $(p,\mu)$-Laplacian is a bounded continuous operator from $W^{1,p}_0(\O)$ to $ W^{-1,p'}(\O)$, and for $p<2$ we have the continuous embedding $W^{1,2}_0(\Omega)\hookrightarrow W^{1.p}_0(\O)$.
 Concerning the \emph{mollified}-convective term, if we consider a
 sequence $\lbrace v_n \rbrace \subset W^{1,2}_0(\Omega)$ such that $v_n \to v \text{ in } W^{1,2}_0(\Omega)$, we have
 $$(J_\mu (v_n) \cdot \nabla v_n, w)  \to  (J_\mu (v) \cdot \nabla v, w) \hspace{1cm } \forall \, w\in W^{1,2}_0(\Omega).$$
 Indeed,
\begin{align*}
 & \left| (J_\mu (v_n) \cdot \nabla v_n, w) - (J_\mu (v) \cdot \nabla v, w)  \right| \\
& = \left| ( (J_\mu (v_n)-J_\mu (v)) \cdot \nabla v_n, w)  +  (J_\mu (v) \cdot \nabla(v_n - v), w)  \right|\\
& \leq || J_\mu (v_n)-J_\mu (v)||_{L^{\infty}(\O)} || \nabla v_n||_{ L^{2}(\Omega)} || w ||_{ L^{2}(\Omega)} \\
& \hspace{0.5cm} +  ||J_\mu (v)||_{L^{\infty}(\O)}  ||\nabla(v_n - v)||_{ L^{2}(\Omega)} ||w||_{ L^{2}(\Omega)}\\
& \leq c(\mu) ||v_n-v||_{L^2(\O)}|| \nabla v_n||_{ L^{2}(\Omega)} || w ||_{ L^{2}(\Omega)} \\
& \hspace{0.5cm}  +  c(\mu) ||v||_{L^2(\O)}||\nabla(v_n - v)||_{ L^{2}(\Omega)} ||w||_{ L^{2}(\Omega)},
\end{align*}
and both terms in the last inequality converge to 0, as $n\to + \infty$.
\item[(C.2)]  Using the $(p,\mu)$-Laplacian growth
    properties, Pettis and Fubini Theorems, and the presence of
    $J_\mu$, it is clear that every term appearing in the definition
    of $A$ belongs to $L^1(\O_T).$ For more details see \cite[Section 6]{BMR}.
  \item[(C.3)]  If we fix $w\in W^{1,2}_0(\Omega)$, with $||w||_{W^{1,2}_0(\Omega)} \leq 1$, since $p<2$, 
 {\small \begin{align*}
 \langle A(v), w\rangle _{W^{1,2}_0(\Omega)} \,  &= \nu  (\nabla v,\nabla w) + ((\mu + |\nabla v|^2)^\frac{p-2}{2}\nabla v,\nabla w) + (J_\mu (v) \cdot \nabla v, w ) \\
& \leq \nu ||\nabla v||_{ L^{2}(\Omega)} + c ||\nabla v ||_{L^{p}(\Omega)}^{p-1}+ ||J_\mu(v)||_{L^{\infty}(\Omega)} ||\nabla v||_{ L^{2}(\Omega)} \\
& \leq ||v||_{W^{1,2}_0(\Omega)}+ c ||\nabla v||^{p-1}_{ L^{2}(\Omega)} + c(\mu) ||v||_{ L^{2}(\Omega)} ||v||_{W^{1,2}_0(\Omega)} \\
&\leq 1\!+\! ||v||_{W^{1,2}_0(\Omega)}\!+\! 1\!+\! ||v||_{W^{1,2}_0(\Omega)} \! +\! c(\mu) ||v||_{ L^{2}(\Omega)} (1\!+\!||v||_{W^{1,2}_0(\Omega)} ) \\
&\leq ( 1+ ||v||_{W^{1,2}_0(\Omega)}) (2 + c(\mu) ||v||_{ L^{2}(\Omega)} ).
   \end{align*}}
So, we choose $\mathcal{B}(||v||_{ L^{2}(\Omega)}) := (2 + c(\mu)||v||_{ L^{2}(\Omega)}),$ $\beta(t) \equiv \alpha (t) \equiv 1,$ $\gamma(t) \equiv 0.$
\item[(C.4)] The property follows immediately as in \cite{BMR}, together with the compactness of the mollified convective term, which can be established in the same way as for the convective term.
\item[(C.5)] We have
\begin{align*}
\langle Av,v\rangle _{W^{1,2}_0(\Omega)} & = \nu  (\nabla v,\nabla v) + ((\mu + |\nabla v|^2)^\frac{p-2}{2}\nabla v , \nabla v)+ (J_\mu (v) \cdot \nabla v, v) \\
&\geq \nu ||\nabla v||^2_{ L^{2}(\Omega)} +   (J_\mu (v) \cdot \nabla v, v)\\
&\geq \nu || v||^2_{W^{1,2}_0(\Omega)} +   (J_\mu (v) \cdot \nabla v, v),
\end{align*}
and for the convective term we get
$$(J_\mu (v) \cdot \nabla v, v) \geq - \left|(J_\mu (v) \cdot \nabla v, v)\right| \geq - ||J_\mu(v)||_{L^{\infty} (\Omega)} || \nabla v||_{ L^{2}(\Omega)} ||v||_{ L^{2}(\Omega)}.$$
By Young inequality, and $||J_\mu(v)||_{L^{\infty}(\O)} \leq c(\mu)
||v||_{ L^{2}(\Omega)}$, we deduce $$\langle Av,v\rangle _{W^{1,2}_0(\Omega)} \geq  \frac{\nu}{2} || v||^2_{W^{1,2}_0(\Omega)}  - c(\nu, \mu) ||v||_{ L^{2}(\Omega)}^4. $$
\end{itemize}
So, in view of Proposition \ref{prop13}, our operator is Bochner
pseudo-monotone.
Thus, Theorem \ref{theo1} yields the assertion. 
\chiu
\subsection{\normalsize The maximum principle for $(\ref{PFepv})$.}
By the duality technique, that is a suitable modification of the one employed,  in \cite{Mar} for a quasi-linear Laplacian system with the same right-hand side as in system (\ref{PF}), we can prove that the  solution of $(\ref{PFepv})$ in Proposition $\ref{existenceL}$ satisfies the maximum principle.
\begin{prop}\label{ulinfty} {\sl 
Let $v^{\mu,\nu}$ be the solution of ($\ref{PFepv})$ on the interval $(0,T^*)$ corresponding to
$v_\circ\in L^{\infty}(\O)$, and $ \mu, \nu >0$. Then \begin{equation}\label{mp}
||v^{\mu, \nu}(t)||_{L^{\infty}(\Omega)} \leq ||v_\circ||_{L^{\infty}(\Omega)}, \hspace{0.5 cm}\text{ for all } \, t \in [0,T^*].
\end{equation}}
\end{prop}
\subsubsection{\normalsize The dual problem}
We introduce an approximating dual problem which will be crucial in
the proof of Proposition \ref{ulinfty}. 
For $v$ given, such that
\be v\in L^{\infty}(0,T;  L^{2}(\Omega)) \cap L^2(0,T; {W^{1,2}_0(\Omega)}), \label{rv}\ee
we denote by $\widetilde J_{\eta}(\nabla v)$ the convolution
$\widetilde \rho_\eta * \nabla v$, where $\widetilde \rho _\eta$ is a
time-spatial mollification kernel. 
For a suitable function $g$, we set $\widehat{g}
  (s):= g(t-s)$.
Moreover, we define, a.e. in $s\in (0,t)$,
\be\label{beeta} (B_\eta(s,x))_{i\beta j\gamma}:=
\delta_{ij}\,\delta_{\beta\gamma}\,\widehat{a}_\eta(\mu, v)(s,x).
\ee
 For $\eta, \mu>0$, and $\nu>0$, we analyse semigroup properties for the following parabolic system 
{\small \begin{align} \partial_s \vp-\nu\Delta \vp- \nabla\cdot
(B_\eta(s,x)\nabla \vp) &= \vp\nabla\cdot J_\mu(\widehat{v})+J_\mu(\widehat{v})\cdot \nabla\vp\,,&&\textrm{ in }(0,t)\times\O, \nonumber \\
\vp(s,x)& =0,&&\textrm{ on }
(0,t)\times\po, \nonumber\\
\vp(0,x)& =\vp_\circ(x),&&\mbox{ on
}\{0\}\times\O, \label{AD1} 
\end{align}}%
where $J_\mu (\widehat{v})$ is a spatial mollifier. The system (\ref{AD1}) with $J_\mu$ replaced by $H_\mu$ is studied in \cite{CDF}\footnote{We take this opportunity to point out a misprint in \cite{CDF}, where in system (2.7) $H_\mu(v)$  has to be replaced by $H_\mu(\widehat{v})$.}.
\begin{lemma}\label{existence}{\sl Assume that  $\eta, \mu>0$, $\nu>0$,  and let $\vp_\circ(x)\in C_0^{\infty}(\O)$. For ${v \in L^\infty (0,T;L^2(\O)) \cap L^2(0,T;W^{1,2}_0(\O))}$, and for all $t\in (0,T)$, there exists
a unique solution $\vp^{\nu, \mu, \eta }$ of $(\ref{AD1})$, such that $\vp^{\nu, \mu, \eta }\in
 L^2(0,t;W^{1,2}_0(\Omega) \cap W^{2,2}(\O))$, $\partial_s \vp^{\nu, \mu, \eta } \in
L^2(0, t;  L^{2}(\Omega))$.
 }\end{lemma}
 \Pr
 The existence of $\vp^{\nu, \mu, \eta }$ is a consequence of \cite[Theorem 2.1]{Den} with $p=2$ and $m=1$. In fact, if we consider
 $\mathcal{A}(s,x,D) = \displaystyle \sum_{|i|\leq 2} a_i (s,x) D^i$, with
 $$a_0(s,x) := \nabla \cdot J_\mu(\widehat{v}), \; \; \;\; a_1(s,x) := J_\mu(\widehat{v}) + \nabla B_\eta, \; \; \; \; a_2(s,x):= \nu + B_\eta;$$
$\mathcal{B}(s,x,D)=I$, where $I$ is the identity tensor, $g_0(s,x) \equiv 0$, and $f\equiv 0$, one readily verifies that the assumptions of this Theorem hold.
\chiu
\begin{lemma}\label{adgt}{\sl Assume that $\nu >0$, $\eta, \mu>0$, and $t>0$. Let   $\vp_\circ(x)\in C_0^{\infty}(\O)$
  and let $\vp^{\nu, \mu, \eta }$ be the unique solution of $(\ref{AD1})$. Then,  for all  $s\in [0,t]$ there holds 
  \be\label{max11a}\|\vp^{\nu, \mu, \eta }(s)\|_{L^1(\O)}
  \leq \|\vp_\circ\|_{L^1(\O)}.
\ee 
 }\end{lemma}
 \Pr
The proof proceeds along the same lines as \cite[Lemma 2.2]{CDF}, with
$J_\mu$ replacing $H_\mu$. For the sake of completeness, we provide
the details. We will write $\vp$ in place of $\vp^{\nu, \mu, \eta }$
for the sake readability.
Let us multiply $(\ref{AD1})_1$ by $\vp (\delta+|\vp|^2)^{-\frac
{1}{2}}$, for some $\delta\in (0,1)$, and integrate by parts on $(0,s)\times \O$, with $s\in [0,t]$. 
Then 
\begin{align}
&\|(\delta+|\vp(s)|^2)^\frac 12\|_{L^1(\O)}+\nu \int_0^s\!\int_\O
(\delta+|\vp|^2)^{-\frac {1}{2}}|\nabla\vp|^2 dx d\tau \nonumber \\ 
&\quad -\nu\!\int_0^s\!\int_\O (\delta+|\vp|^2)^{-\frac {3}{2}}(\nabla\vp\cdot
\vp)^2 dx d\tau \dy + \int_0^s\!\int_\O B_\eta(\tau, x)\,(\delta+|\vp|^2)^{-\frac
{1}{2}}\,|\nabla\vp|^2 dx d\tau
\nonumber \\ 
&\quad  -\int_0^s\!\int_\O B_\eta(\tau, x)\,(\delta+|\vp|^2)^{-\frac {3}{2}}\,(\nabla\vp\cdot \vp)^2 dx d\tau \nonumber \\
 &= \|(\delta+|\vp(0)|^2)^\frac 12\|_{L^1(\O)} \nonumber\\
&\quad + \int_0^s\!\int_\O\bigg[ (\delta+| \vp|^2)^{-\frac {1}{2}} |\vp|^2\nabla\cdot J_\mu(\widehat{v})+  J_\mu(\widehat{v})\cdot \nabla (\delta+| \vp|^2)^\frac {1}{2} dx d\tau\bigg].\label{max13}
  \end{align}
  Observing that 
  $$(\nu+B_\eta(\tau, x))\bigg[(\delta+|\vp|^2)^{-\frac
{1}{2}}\,|\nabla\vp|^2-(\delta+|\vp|^2)^{-\frac {3}{2}}\,(\nabla\vp\cdot \vp)^2\bigg] \geq 0,$$
and  integrating by parts the last term on the right-hand side of $(\ref{max13})$, we find
\begin{align*}
  &\|(\delta+|\vp(s)|^2)^\frac 12\|_{L^1(\O)}\\
  &\leq \dy  
 \|(\delta+|\vp_0|^2)^\frac 12\|_{L^1(\O)}\\ 
&\quad  + \int_0^s\!\int_\O  \bigg[ (\delta+| \vp|^2)^{-\frac {1}{2}} |\vp|^2\nabla\cdot J_\mu(\widehat{v})- \nabla\cdot J_\mu(\widehat{v}) (\delta+| \vp|^2)^\frac {1}{2} \bigg]dxd\tau.
  \end{align*}
Since  
$$\lim_{\delta\to 0}\bigg[ (\delta+| \vp|^2)^{-\frac {1}{2}} |\vp|^2\nabla\cdot J_\mu(\widehat{v})- \nabla\cdot J_\mu(\widehat{v}) (\delta+| \vp|^2)^\frac {1}{2} \bigg]=0,$$
and 
\begin{align*}
&\bigg| (\delta+| \vp|^2)^{-\frac {1}{2}} |\vp|^2\nabla\cdot J_\mu(\widehat{v}) - \nabla\cdot J_\mu(\widehat{v}) (\delta+| \vp|^2)^\frac {1}{2} \bigg|  \\
&= \bigg|\frac{|\vp|^2}{ (\delta+| \vp|^2)^{\frac {1}{2}}} -(\delta+| \vp|^2)^\frac {1}{2} \bigg| |\nabla\cdot J_\mu(\widehat{v})| 
 \leq |\nabla\cdot J_\mu(\widehat{v})|,  
\end{align*}
for $\delta$ sufficiently small, we can apply the Lebesgue dominated convergence theorem, and get $(\ref{max11a})$.\chiu
\subsection{\normalsize Proof of Proposition \ref{ulinfty}}
\Pr Let $t\in [0,T^*]$ be arbitrary but fixed, and $\nu, \mu, \eta>0$ fixed. Writing $v$ in place of $v^{\mu, \nu}$, we  consider  the solution $\vp^\eta(s,x)$
of system $(\ref{AD1})$ corresponding to this $v$ and an initial datum $\vp_\circ\in
C_0^\infty(\O)$, where $B_\eta$ is given by $(\ref{beeta})$ and
$J_\mu(\widehat{v})$ the spatial mollifier of the time shifted $v$.  
Then, using $\widehat\vp^\eta(\tau)$ as test
function in the weak formulation of $(\ref{PFepv})$, we
have
 \begin{align} 
 &(v(t),\vp_\circ) -(v_\circ,\vp^\eta(t))
-\int_0^t (v(\tau), \partial_\tau \widehat\vp^\eta(\tau))d\tau  =-\, \nu\int_0^t(\nabla v(\tau),\nabla
\widehat\vp^\eta(\tau))d\tau \nonumber
\\
& -\int_0^t(a(\mu,v(\tau)) \nabla
v(\tau),\nabla \widehat\vp^\eta(\tau))d\tau  -\int_0^t (J_\mu(v)\cdot \nabla v, \widehat\vp^\eta(\tau))d\tau. \label{a8}
\end{align}
We write the second 
term on the right-hand side of $(\ref{a8})$, as
\begin{align}&\int_0^t(a(\mu,v(\tau))
\nabla v(\tau),\nabla \widehat\vp^\eta(\tau))d\tau=\!\!\! \int_0^t(
\nabla v(\tau),a_\eta(\mu,v(\tau))\nabla
\widehat\vp^\eta(\tau))d\tau \nonumber\\
&\dy +\int_0^t( \nabla v(\tau),\nabla
\widehat\vp^\eta(\tau))[a(\mu,v(\tau)) -a_\eta(\mu, v(\tau))]\,d\tau.\label{mts2}\end{align}
Further, we integrate by parts the third term on the right-hand side of $(\ref{a8})$, and we get that
\be\label{mts1}\ba{ll}\dy \int_0^t\!\! (J_\mu(v)\cdot\! \nabla v, \widehat\vp^\eta(\tau))d\tau=-\!
\!\int_0^t\! (\nabla\!\cdot J_\mu(v) v, \widehat\vp^\eta(\tau))d\tau-\int_0^t\! (v,J_\mu(v)\!\cdot\nabla\widehat\vp^\eta(\tau))d\tau.
\ea\ee
Since $\vp^\eta$ is the solution of $(\ref{AD1})$, there holds
\begin{align}
   & \int_0^t (\partial_s \vp^\eta(s),  \hat{v}(s)) + \nu(\nabla \vp^\eta (s), \nabla \hat{v}(s)) + (B_\eta(s,x) \nabla \vp^\eta (s), \nabla \hat{v}(s))\, ds  \notag\\
   & = \int_0^t (\vp^\eta (s) \nabla \cdot J_\mu (\hat{v})(s), \hat{v}(s)) + (J_\mu(\hat{v})(s)\cdot \nabla \vp^\eta(s), \hat{v}(s))\, ds; \label{weakform}
\end{align}
using the transformation $s=t-\tau$ in \eqref{weakform}, and taking into account $(\ref{mts2})$ and $(\ref{mts1})$,  identity
$(\ref{a8})$ becomes 
\begin{align}
 (v(t),\vp_\circ)&=(v_\circ,\vp^\eta(t)) + \int_0^t( \nabla
v(\tau),\nabla \widehat\vp^\eta(\tau))[a(\mu, v(\tau))-a_\eta(\mu,
v(\tau))]\,d\tau \nonumber \\
& 
=:(v_\circ,\vp^\eta(t))+I_\eta. \label{a9}\end{align}%
Using Lemma \ref{comp} with $\nabla h^\eta=\nabla \widehat\vp^\eta\in
L^2(0,t; L^{2}(\Omega))$ and  $\nabla\psi= \nabla v\in L^{2}(0,t;  L^{2}(\Omega))$,
for $\nu>0$ fixed, 
the integral
$I_\eta$ goes to zero, as $\eta$ goes to zero, along a subsequence.
Finally, using $(\ref{max11a})$ and then passing to the limit as $\eta$
tends to zero  in $(\ref{a9})$, along a subsequence, we get that
$$|(v(t),\vp_\circ)| \leq \|v_\circ\|_{L^\infty(\O)}\,\|\vp_\circ\|_{L^1(\O)}, \
\text{ for all } \vp_\circ\in C_0^\infty(\O).$$
A density argument extends this inequality to any $\varphi_\circ\in
L^1(\O)$.  This in turn implies 
$$\|v(t)\|_{L^\infty(\O)}=\!\!\sup_{\vf_\circ\in L^1(\O) \atop
||\vf_\circ||_{L^1(\O)}=1}|(v(t),\vf_\circ)| \leq
\|v_\circ\|_{L^\infty(\O)}, \text{ for all } t\in [0,T^*].$$
\chiu
\subsection{\normalsize Global existence of weak solution of $(\ref{PFepv})$.}
\begin{prop}
\label{existencemunu}
{\sl Let be $\nu, \mu>0$. Assume that $v_\circ$
belongs to $L^{\infty}(\O)$. Then, for any $T>0$ there exists a global solution $v^{\mu,\nu}$ of system $(\ref{PFepv})$ satisfying
\begin{equation}
||v^{\mu, \nu}(t)||_{L^{\infty}(\Omega)} \leq ||v_\circ||_{L^{\infty}(\Omega)}, \text{ for all } t\in [0,T];
\label{mmm}
\end{equation}
\begin{equation}\label{ue}
||v^{ \mu, \nu}||^2_{C ([0,T];  L^{2}(\Omega))} + \nu ||\nabla v^{\mu, \nu}||^2_{L^2 (\Omega_T)} + ||\nabla v^{ \mu, \nu}||^p_{L^{p} (\Omega_T)} \leq K_1(||v_\circ ||_{L^{\infty} (\Omega)},T),
\end{equation}
with $K_1(||v_\circ||_{L^{\infty} (\Omega)},T)$ independent of $ \mu, \nu$.}
\end{prop}
\Pr
Due to Proposition\,\ref{existenceL},  there exists a solution $v^{\mu, \nu}$,
 of $(\ref{PFepv})$ in the time interval $(0, T^* )$ with $T^\star = \frac{\nu\mu^{\frac 32} }{2 \|v_\circ\|_2^2}$ corresponding to the initial datum
$v_\circ$. We will write $v$ in place of $v^{\mu, \nu}$ for the sake
of readability. As, by assumption,  $v_\circ\in L^{\infty}(\O)$, from
Proposition \ref{ulinfty} we also know that $v\in
L^{\infty}(\O_{T^*})$ and, that $v$ satisfies for any $t\in[0,T^*]$  the estimate
\be\label{v2i}
\|v(t)\|_{L^\infty(\O)}\leq \|v_\circ\|_{L^\infty(\O)}.\ee
The validity of the above inequality ensures that the solution $v$ actually exists for any $T>0$, since we can apply Proposition \ref{existenceL} with initial time $T^*$ and repeat this process finitely many times. Moreover, Proposition $\ref{ulinfty}$ ensures the validity of $(\ref{mmm})$. \\
 By Remark \ref{rmk1}, we test the weak formulation of $\eqref{weakmunualpha}$ by $v$ and use the  integration by parts formula (cf. $\eqref{intbyp}$), to get that
 \begin{align*}
& ||v(t)  ||^2_{ L^{2}(\Omega)} +  2 \nu \int_0^t ||\nabla v(\tau)||_{ L^{2}(\Omega)}^2 \, d\tau +  2 \int_0^t ||(\mu + |\nabla v(\tau)|^2)^\frac{p-2}{4} \nabla v(\tau)||_{ L^{2}(\Omega)}^2\, d\tau \\
&  =||v_\circ||_{ L^{2}(\Omega)}^2 - \int_0^t (J_\mu (v) \cdot \nabla v,v) \, d\tau.
  \end{align*}
 Moreover, by Hölder's inequality, interpolation, ($\ref{M2}$),  and
 ($\ref{mmm}$) we have
 \begin{align*}
\left|  (J_\mu (v) \cdot \nabla v,v)  \right| & \leq ||J_\mu (v) ||_{L^{\infty} (\Omega)} ||\nabla v||_{L^{p}(\Omega)} || v||_{L^{p'}(\Omega)}\\
&\leq \frac{1}{p} ||\nabla v||_{L^{p}(\Omega)}^p + \frac{1}{p'} ||v_\circ||_{L^{\infty}(\Omega)}^{\frac{2}{p-1}} ||v||^2_{ L^{2}(\Omega)},    
 \end{align*}
 and by using, on the left-hand side, the estimate 
\begin{align}
\int_{\O}|\nabla v|^p\, dx  &=\int_{|\nabla v|^2> \mu}|\nabla v|^p\, dx+\int_{|\nabla v|^2 \leq \mu}|\nabla v |^p\, dx \nonumber \\ 
& \leq 2^\frac{2-p}{2}\!\int_{\O}(\mu + |\nabla v|^2)^\frac{p-2}{2}  \,|\nabla v|^2 dx+\int_{|\nabla v|^2\leq \mu} \mu^\frac p2\, dx \nonumber \\ 
& \leq  2^\frac{2-p}{2} \int_{\O} (\mu + |\nabla v|^2)^\frac{p-2}{2}\, |\nabla v|^2dx+ \mu^\frac p2 \, |\O|, \label{gp}
\end{align}
 we get
 \begin{align*}
& ||v(t)||^2_{ L^{2}(\Omega)} + 2 \nu \int_0^t ||\nabla v(\tau)||_{ L^{2}(\Omega)}^2 \, d\tau + c \int_0^t ||\nabla v(\tau)||_{L^{p}(\Omega)}^p \, d\tau \\
& \leq ||v_\circ||^2_{ L^{2}(\Omega)} + \mu^{\frac{p}{2}} |\Omega| T + \frac{1}{p} ||v_\circ||_{L^{\infty}(\Omega)}^\frac{2}{p-1} \int_0^T ||v(\tau)||_{ L^{2}(\Omega)}^2 \, d\tau.
  \end{align*}
Ultimately, we have
$$ ||v^{ \mu, \nu}||^2_{L^{\infty} (0,T;  L^{2}(\Omega))} + \nu ||\nabla v^{ \mu, \nu}||^2_{L^2 (\Omega_T)} + ||\nabla v^{ \mu, \nu}||^p_{L^{p} (\Omega_T)} \leq K_1(||v_\circ ||_{L^{\infty} (\Omega)}, T),$$
where $K_1$ does not depend on  $ \mu, \nu$. 
{Note that Proposition \ref{contrap} implies \linebreak ${v^{\mu, \nu} \in C([0,T];L^2(\O)).}$}
\chiu
\section{\normalsize Limiting process as $\nu \to 0^+$} \renewcommand{\theequation}{4.\arabic{equation}}
\renewcommand{\thetho}{4.\arabic{tho}}
\renewcommand{\thedefi}{4.\arabic{defi}}
\renewcommand{\therem}{4.\arabic{rem}}
\renewcommand{\theprop}{4.\arabic{prop}}
\renewcommand{\thelemma}{4.\arabic{lemma}}
\renewcommand{\thecoro}{4.\arabic{coro}}
\setcounter{tho}{0}
In this section, we prove that a suitable sequence of solutions to
system $(\ref{PFepv})$ converges, as $\nu \to 0^+$, to a solution of
the following system $(\ref{PFep})$ in the sense of Definition $\ref{wmualpha}$:
\begin{align} \partial_t v-\nabla\cdot\left((\mu+|\nabla
v|^2)^\frac{p-2}{2}\nabla v\right)&= -J_\mu(v)\cdot \nabla v\,,&&\textrm{ in
}(0,T)\times\O, \nonumber \\ 
v(t,x)&=0\,,& &\textrm{ on }
(0,T)\times\po, \nonumber \\
v(0,x)&=v_\circ(x),&&\mbox{ on
}\{0\}\times\O.\label{PFep}
\end{align} 
\begin{defi}\label{wmualpha} {\rm Let $\mu> 0$, $T>0$, and let
    $v_{\circ}\in L^{\infty}(\O)$ be given. A vector field $v\!:(0,T)\times \O\to\R^3$
 is said to be a global weak solution of
system {\rm $(\ref{PFep})$} on $(0,T)$, if
\begin{itemize}
\item[i)] 
$  v\in  L^{p}(0,T; W_0^{1,p}(\O))\cap L^{\infty}(\O_T)$, \par $\partial_t v \in L^{p'}(0,T;W^{-1,{p'}}(\O)) + L^{p}(0,T;L^{p}(\O));$
\item[ii)] $v$ satisfies the integral identity
$$\int_0^t \left[(v,\partial_\tau \psi) -\left(a(\mu,
v)\,\nabla v,\nabla \psi\right)-(J_\mu(v)\cdot \nabla v, \psi)\right]d\tau=(v(t),
\psi(t))-(v_\circ, \psi(0)),$$
for all  $t\in [0,T]$, and all $\psi\in C^\infty([0,T]; C_0^\infty( \O))$; 
\item[iii)]$\dy \lim_{t\to
0^+}\|v(t)-v_\circ\|_{ L^{2}(\Omega)}=0\,.$\end{itemize}}
\end{defi}
\begin{rem} By a density argument and the regularity of solution $v$ in Definition~\ref{wmualpha}, we can consider test functions $\psi$ in ii) such that $\psi\in  L^{p}(0,T; W_0^{1,p}(\O))\cap L^{\infty}(\O_T)$, and $\partial_t \psi \in L^{p'}(0,T;W^{-1,{p'}}(\O)) + L^{p}(0,T;L^{p}(\O))$. \label{rmk2}
\end{rem}
\begin{prop}\label{existencemualpha} 
{\sl Let $\mu>0$. Assume that $v_\circ$
belongs to $L^{\infty}(\O)$. Then, for any $T>0$ there exists a global weak solution $v^\mu$
of system $(\ref{PFep})$ in the sense of Definition\;$\ref{wmualpha}$.
In particular, 
\begin{equation}\label{ues}
\dy \|v^\mu\|^2_{C([0,T]; L^{2}(\Omega))}+\|\nabla v^\mu\|^p_{L^{p}(0,T; L^{p}(\O))}\leq K_1(\|v_\circ\|_\infty, T)\,. \end{equation}
Moreover, $v$ satisfies the maximum principle
\begin{equation}\label{mp2}
||v^\mu(t)||_{L^{\infty}(\Omega)} \leq ||v_\circ||_{L^{\infty} (\Omega)}, \;\; \; \text{for all }  t \in [0,T].
\end{equation}}
\end{prop}
\Pr 
Let us fix $\mu >0$. Then, for any fixed but arbitrary $\nu>0$ there
exists a solution $v^{\mu, \nu}$ of $(\ref{PFepv})$ as in Proposition
$\ref{existencemunu}$  corresponding to $v_\circ \in
L^{\infty}(\O)$. We will write $v^\nu$ in place of $v^{\mu, \nu}$ for
the sake of readability.\\
 \noindent{\bf Convergence property} - Since $v^\nu$ is a solution of the system $(\ref{PFepv})$ in the sense of Definition $\ref{weakmunualpha}$ the integration by parts formula \eqref{intbyp} yields 
 \begin{align}
 \int_0^t \langle \partial_\tau v^\nu, \psi \rangle_{W^{1,2}_0(\O)} \,d\tau&=  \,- \nu\int_0^t (\nabla v^\nu, \nabla \psi) + (a(\mu, v^\nu) \nabla v^\nu, \nabla \psi) \, d\tau
\nonumber \\
&\quad - \int_0^t  (J_\mu(v^\nu)\cdot \nabla v^\nu, \psi)\, d\tau , \label{weakf} 
 \end{align}
for all $\psi\in C^\infty([0,T]; C_0^\infty( \O))$. 
Note that for fixed $\mu>0$, the families
\begin{align*}
\nu\Delta v^\nu & \in L^{p}(0,T;W^{-1,p} (\Omega)),\\
\nabla\cdot\left((\mu+|\nabla
v^\nu|^2)^\frac{p-2}{2}\nabla v^\nu\right) & \in L^{p'}(0,T;W^{-1,p'} (\Omega)), \\
J_\mu(v^\nu)\cdot \nabla v^\nu & \in L^{p}(0,T;L^{p}(\Omega))
\end{align*}
are uniformly bounded with respect to $\nu>0$ due to $(\ref{ue})$. 
Since $p<2$, we have $  L^{p'}(0,T;W^{-1,p'} (\Omega)) \subset
L^{p}(0,T;W^{-1,p} (\Omega))$,  $L^{p}(0,T;L^{p} (\Omega)) \subset
L^{p}(0,T;W^{-1,p} (\Omega))$. Thus, we get that
$\partial_t v^\nu $ is uniformly bounded in
$L^{p}(0,T;W^{-1,p}(\Omega))$ with respect to $\nu>0$. Since
$W^{1,p}_0(\Omega) \hookrightarrow \hookrightarrow
L^{p}(\Omega)\hookrightarrow  W^{-1,p}(\Omega)$ the Aubin-Lions
Theorem (cf.~Theorem \ref{aubin}) implies that there exists a
subsequence and a $v=v^\mu \in L^{p}(0,T; L^{p}(\Omega))$ such  that 
\begin{equation}\label{al}
v^\nu \to v \text{ in }L^{p}(0,T; L^{p}(\Omega)).
\end{equation}
Summarising, from the properties given in $(\ref{mmm})$, $(\ref{ue})$ and $(\ref{al})$, also using interpolation, we can choose a subsequence $v^\nu$ such that
 \begin{align}
 \label{A}  v^\nu &\to v & &\text{ in } L^q (0,T; L^q (\Omega)) \text{ strongly, for all }  q  \in [1, +\infty);\\
\label{B} v^\nu &\overset{\ast}{\rightharpoonup} v & &\text{ in } L^{\infty} (\Omega_T) \text{ weakly-*};\\
\label{BB} v^\nu &\overset{\ast}{\rightharpoonup} v & &\text{ in } L^{\infty}(0,T;L^2 (\Omega)) \text{ weakly-*};\\
\label{C} v^\nu &\rightharpoonup v & &\text{ in } L^{p}( 0,T;W^{1,p}_0(\Omega)) \text{ weakly};\\
\label{D}a(\mu, v^\nu) \nabla v^\nu &\rightharpoonup\chi & &\text{ in } L^{p'} (0,T; L^{p'}(\Omega))\text{ weakly};\\
\label{E} v^\nu_t&\rightharpoonup \partial_t v & &\text{ in } L^{p'} (0,T; W^{-1,p'}(\Omega)) + L^p(0,T;L^p(\O))  \text{ weakly}.
\end{align}
Note that \eqref{B}, \eqref{C}, \eqref{E} and Proposition \ref{con}
imply that
\begin{align}\label{c2}
  v \in C([0,T];L^2(\O))\,.
\end{align}
\noindent{\bf Limiting process} -  If we use $\psi \in C^\infty_0(\Omega_T)$ as test function in  the weak formulation of $(\ref{PFepv})$ in Definition $\ref{weakmunualpha}$, we have
 \be \label{weakff}\int_0^T  (v^\nu, \partial_\tau \psi) - \nu (\nabla v^\nu,\nabla \psi) -(a(\mu, v^\nu) \nabla v^\nu, \nabla \psi)  \, d\tau \; =  \int_0^T (J_\mu(v^\nu)\cdot v^\nu, \psi) \, d\tau.\ee
 Now,  as $\nu \to 0^+$, by $(\ref{A})$
 $$\int_0^T (v^\nu, \partial_\tau \psi) \, d \tau \to \int_0^T (v, \partial_\tau \psi) \, d \tau; $$
 by $(\ref{ue})$
 $$\nu\int_0^T (\nabla v^\nu,\nabla \psi) \, d \tau \to 0,$$
indeed
\begin{align*}
    \left| \nu\int_0^T (\nabla v^\nu,\nabla \psi) \, d \tau \right| & \leq \nu ||\nabla v^\nu||_{L^{2}(\Omega_T)} ||\nabla \psi||_{L^{2}(\Omega_T)}\\
    &\leq \nu^{\frac 12} K_1^{\frac 12} ||\nabla \psi||_{L^{2}(\Omega_T)} \to 0,  \text{ as } \nu \to 0^+;
\end{align*}
  by $(\ref{D})$
  $$\int_0^T(a(\mu, v^\nu) \nabla v^\nu, \nabla \psi)   \, d\tau \to \int_0^T(\chi, \nabla \psi)   \, d\tau;$$
  by $(\ref{A})$, $(\ref{C})$, and $(\ref{M1})$
  $$\int_0^T (J_\mu(v^\nu)\cdot \nabla v^\nu, \psi) \, d\tau \to \int_0^T (J_\mu(v)\cdot  \nabla v, \psi) \, d\tau.$$
So, we get $v\in L^{p}(0,T;W^{1,p}_0(\Omega)) \cap L^{\infty}(\Omega_T)$, with
\begin{equation}
\partial_\tau v =  \nabla \cdot \chi - J_\mu(v)\cdot \nabla v \in L^{p'}(0,T;W^{-1,p'}(\Omega)) + L^{p}(0,T;L^{p}(\Omega)),
\label{vt}
\end{equation} and we need to show that
  $\chi = a(\mu, v)\nabla v.$ To this end we define
 \be X_\nu := \int_0^T (a(\mu, v^\nu)\nabla v^\nu - a(\mu, \varphi)\nabla \varphi, \nabla v^\nu - \nabla \varphi) \, d \tau. \label{chimon}\ee
  From the properties of the $(p,\mu)$-Laplacian it is well known that $X_\nu \geq 0$,  for all $ \varphi \in L^{p}(0,T;W^{1,p}_0(\Omega)).$
 Note that from $\eqref{weakf}$ follows that, for
  each $\nu$, $\partial_t v^\nu \in L^2(0,T;W^{-1,2}(\O)),$ where we
  used $a(\mu, v)\leq \mu^{\frac{p-2}{2}}$, the properties of
  mollifier, and \eqref{ue}. Thus, due to Remark \ref{rmk1}, and  the
  integration by parts formula (cf.~\eqref{intbyp}) we get:
\begin{align} 
&\int_0^T (a(\mu, v^\nu)\nabla v^\nu, \nabla v^\nu) \, d\tau \,   =  \, -\frac{1}{2} ||v^\nu (T)||_{ L^{2}(\Omega)}^2 + \frac{1}{2} ||v_\circ||_{ L^{2}(\Omega)}^2 \nonumber \\
& - \nu \int_0^T(\nabla v^\nu, \nabla v^\nu) \, d\tau - \int_0^T (J_\mu (v^\nu) \cdot \nabla v^\nu, v^\nu)\, d\tau. \label{lsnu}
\end{align}
 \\ We have
\begin{equation}\label{lsl}
\limsup_{\nu \to 0^+} \left(- \nu \int_0^T(\nabla v^\nu, \nabla v^\nu) \, d\tau \right) \leq 0,
\end{equation} 
and, by $(\ref{A})$, $(\ref{C})$, and $(\ref{M1})$
\begin{equation}\label{lsc}
\limsup_{\nu \to 0^+} \left(\int_0^T (J_\mu (v^\nu) \cdot \nabla v^\nu, v^\nu)\, d\tau \right) = \int_0^T (J_\mu (v) \cdot \nabla v, v)\, d\tau.
\end{equation} 
Moreover, we prove that it is possible to choose the subsequence $v^\nu$ in such a way
$v^\nu(T)\rightharpoonup v(T) \text{ in }  L^{2}(\Omega)$, and
consequently we get
\begin{equation}\label{lsT}\liminf_{\nu \to 0^+} ||v^\nu(T)||_{ L^{2}(\Omega)} \geq ||v(T)||_{ L^{2}(\Omega)}.
\end{equation}
To this end we note that $v^\nu \in C(\left[ 0,T\right];  L^{2}(\Omega))$, and that
$v^\nu(T)$ is uniformly bounded in $ L^{2}(\Omega)$ (cf.~\eqref{mmm}).
Thus, there exists a subsequence, $v^\nu (T)$, such that 
$v^\nu(T)\rightharpoonup \xi_T$ in $  L^{2}(\Omega).$
Using $\psi (t,x) = \phi(t) \omega(x) \in C^\infty([0,T];
C^\infty_0(\Omega)),$ with $\phi (0)=1, \phi (T)=0$ or $\phi (0)=0, \phi (T)=1$, as test function in \eqref{weakff}, and passing to the limit w.r.t. $\nu \to 0^+$, we obtain 
\begin{align*}
&\int_0^T (v, \partial_\tau \psi) \, d\tau - \int_0^T(\chi, \nabla \psi) \, d\tau \\ &= \phi(T) (\xi_T, \omega(x)) - \phi(0)(v_\circ, \omega(x)) + \int_0^T( J_\mu (v) \cdot \nabla v, \psi) \, d\tau.
\end{align*}
On the other hand, using $(\ref{vt})$ together with Proposition $\ref{ibp}$, we get 
\begin{align*}
&\int_0^T (v, \partial_\tau \psi) \, d\tau - \int_0^T(\chi, \nabla \psi) \, d\tau  \\&= \phi(T) (v(T), \omega(x)) - \phi(0)(v(0), \omega(x)) + \int_0^T( J_\mu (v) \cdot \nabla v, \psi) \, d\tau.
\end{align*}
So, we deduce that $\xi_T=v(T)$ in $ L^{2}(\Omega)$, and we get
$(\ref{lsT})$. Moreover, we have $v_\circ = v(0)$ in $ L^{2}(\Omega).$

Coming back to $(\ref{lsnu})$, using $(\ref{lsl})$, $(\ref{lsc})$, $(\ref{lsT})$, we have
\begin{align*}
&\limsup_{\nu \to 0^+}\left( \int_0^T (a(\mu, v^\nu)\nabla v^\nu,
  \nabla v^\nu) \, d\tau \right)
\\
&\leq -\frac{1}{2}||v(T)||_{ L^{2}(\Omega)}^2 + \frac{1}{2}||v_\circ||_{ L^{2}(\Omega)}^2 - \int_0^T(J_\mu(v)\cdot \nabla v, v)\, d\tau.
\end{align*}
Moreover, from \eqref{chimon} follows
\begin{align*}
 0\leq \limsup_{\nu \to 0^+} X_\nu &\leq  -\frac{1}{2}||v(T)||_{ L^{2}(\Omega)}^2 + \frac{1}{2}||v_\circ||_{ L^{2}(\Omega)}^2 - \int_0^T(J_\mu(v)\cdot \nabla v, v)\, d\tau\\
 &\quad -\int_0^T(\chi, \nabla \varphi)\,d\tau -\int_0^T ( a(\mu, \varphi)\nabla \varphi, \nabla v - \nabla \varphi) \, d \tau,
 \end{align*}
and, since we have proved that $v(0)=v_\circ$, Proposition $\ref{ibp}$ guaranties 
$$-\frac{1}{2}||v(T)||_{ L^{2}(\Omega)}^2 + \frac{1}{2}||v_\circ ||_{ L^{2}(\Omega)}^2 \; = -\int_0^T \langle\partial_\tau v , v\rangle _{W^{1,p}_0(\Omega) \cap L^{p'}(\Omega)} \, d\tau. $$
This implies
$$0 \leq \int_0^T (\chi - a(\mu, \varphi)\nabla \varphi, \nabla v - \nabla \varphi) \, d\tau, \; \;\text{for all } \varphi \in L^{p}(0,T; W^{1,p}_0(\Omega)),$$
so, by Minty Trick, $\chi = a(\mu,v)\nabla v.$\\
\textbf{Maximum principle} - 
In view of (\ref{mmm}) there exists for each fixed $t \in [0,T]$ a
subsequence $v^{\nu '}(t)$ of $v^{\nu }(t)$ such that $v^{\nu
  '}(t)\overset{\ast}{\rightharpoonup} \xi_t$ in
$L^{\infty}(\Omega)$. Moreover, the uniform boundedness of $v^\nu (t)$
in $C([0,T]; L^{2}(\Omega))$ (cf.~\eqref{ue}) yields that there exists
a subsequence $v^{\nu ''}(t)$ of $v^{\nu '}(t)$ such that $v^{\nu
  ''}(t)\rightharpoonup \tilde{\xi}_t$ in $ L^{2}(\Omega)$. By
uniqueness of the weak limits we have $\tilde{\xi}_t= \xi_t$. In view
of \eqref{c2} we also have that ${v \in C([0,T]; L^{2}(\Omega))}$, and thus,
we want to show that $\xi_t = v(t)$ in $L^2(\O)$.
To this end we consider $\tilde{\psi}(\tau,x)= \phi(\tau) \omega(x)
\in C^\infty ([0,T];C^\infty_0(\Omega))$, with $\phi(t)=1$ and
$\phi(0)=0$, and define $\psi(\tau,x) := \tilde{\psi}(\tau,x)
\chi_{\left[0,t\right]} (\tau)$. Choosing $\psi=v^{\nu '}(t)$ in the
weak formulation of Definition $\ref{weakmunualpha}$ for $\nu'$, and passing to the limit, we get
$$\int_0^t (v,\partial_\tau \tilde{\psi}) \, d\tau -\int_0^t ((\mu + |\nabla v|^2)^{\frac{p-2}{2}} \nabla v, \nabla \tilde{\psi}) \, d\tau\, = \int_0^t (v\cdot \nabla v, \tilde{\psi}) \, d\tau - (\xi_t, w).$$
Using $(\ref{vt})$ and Proposition $\ref{ibp}$  we get 
$$\int_0^t (v,\partial_\tau \tilde{\psi}) \, d\tau -\int_0^t ((\mu + |\nabla v|^2)^{\frac{p-2}{2}} \nabla v, \nabla \tilde{\psi})) \, d\tau\, = \int_0^t (v\cdot \nabla v, \tilde{\psi}) \, d\tau - (v(t), w),$$
which implies $\xi_t = v(t)$ in $L^2(\O)$ and
$||v(t)||_{L^{\infty}(\Omega)} \leq
||v_\circ||_{L^{\infty}(\Omega)}$, by $(\ref{mmm})$ and
semi-continuity of $L^{\infty}$-norm. Using the same argument for each $t \in \left[0,T\right]$, we get ($\ref{mp2}$).
\\\textbf{Initial datum} - Condition $iii)$ is satisfied. Indeed, due
to \eqref{c2} we have ${v \in C([0,T]; L^{2}(\Omega))}$ and we have proved that $v(0)=v_\circ$.
\chiu
\section{\normalsize Proof of Theorem $\ref{exis}$ }
\setcounter{tho}{0}
\renewcommand{\theequation}{5.\arabic{equation}}
\renewcommand{\thetho}{5.\arabic{tho}}
\renewcommand{\thedefi}{5.\arabic{defi}}
\renewcommand{\therem}{5.\arabic{rem}}
\renewcommand{\theprop}{5.\arabic{prop}}
\renewcommand{\thelemma}{5.\arabic{lemma}}
\renewcommand{\thecoro}{5.\arabic{coro}}
In this section, we obtain a solution of $(\ref{PF})$, for $v_\circ \in L^{\infty}(\O)$, as the limit of a sequence of solutions for $(\ref{PFep})$, as $\mu \to 0^+$.\\
\Pr For any fixed $\mu>0$, denote by $v^\mu$  a solution of the system
$(\ref{PFep})$ as in Proposition $\ref{existencemualpha}$ corresponding to $v_\circ \in L^{\infty}(\O)$.\\
 \noindent{\bf Convergence property} - Since $v^\mu$ is a solution of
 the system $(\ref{PFep}) $ as in the sense of Definition
 $\ref{wmualpha}$ the integration by parts formula in the Proposition $\ref{ibp}$,  yields
\be
\label{weakf1}
\int_0^t \langle \partial_\tau v^\mu, \psi \rangle_{W^{1,p}_0(\Omega) \cap L^{p'}(\Omega)} \,d\tau= - \int_0^t  (a(\mu, v^\mu) \nabla v^\mu, \nabla \psi) + (J_\mu(v^\mu)\cdot \nabla v^\mu, \psi)\, d\tau  
 \ee
for all $\psi\in C^\infty([0,T]; C_0^\infty( \O))$.
Note that the families
\begin{align*}
\nabla\cdot\left(a(\mu, v^\mu)\nabla v^\mu\right) & \in L^{p'}(0,T;W^{-1,p'} (\Omega)), \\
J_\mu(v^\mu)\cdot \nabla v^\mu & \in L^{p}(0,T;L^{p}(\Omega))
\end{align*}
are uniformly bounded with respect to $\nu>0$ due to $(\ref{ues})$. 
Since $p<2$, we have $  L^{p'}(0,T;W^{-1,p'} (\Omega)) \subset
L^{p}(0,T;W^{-1,p} (\Omega))$, $L^{p}(0,T;L^{p} (\Omega)) \subset
L^{p}(0,T;W^{-1,p} (\Omega))$. Thus, we get that $\partial_t v^\mu$
is uniformly bounded in $L^{p}(0,T;W^{-1,p}(\Omega))$ with respect to $\nu>0$. Since
$W^{1,p}_0(\Omega) \hookrightarrow \hookrightarrow L^{p}(\Omega)\hookrightarrow  W^{-1,p}(\Omega)$ 
the Aubin-Lions Theorem (cf.~Theorem \ref{aubin}) implies that there
exists a subsequence and a $v \in L^{p}(0,T; L^{p}(\Omega))$ such that 
\begin{equation}\label{almu}
v^\mu \to v \text{ in }L^{p}(0,T; L^{p}(\Omega)).
\end{equation}
Summarising, from the properties given in $(\ref{ues})$, $(\ref{mp2})$
and $(\ref{almu})$ and $(\ref{al})$, also using interpolation, we can choose a
subsequence $v^\nu$ such that
 \begin{align}
 \label{Amu} v^\mu \to v & \text{ in } L^q (0,T;L^q(\Omega)) \text{ strongly, for all }  q  \in [1, +\infty);\\
\label{Bmu} v^\mu \overset{\ast}{\rightharpoonup} v & \text{ in } L^{\infty} (\Omega_T) \text{ weakly-*};\\
\label{BBmu} v^\mu \overset{\ast}{\rightharpoonup} v & \text{ in } L^{\infty} (0,T; L^2(\Omega) \text{ weakly-*};\\
\label{Cmu} v^\mu \rightharpoonup v & \text{ in } L^{p} (0,T;W^{1,p}_0(\Omega)) \text{ weakly}; \\
\label{Dmu} a(\mu, v^\mu)\nabla v^\mu \rightharpoonup\chi & \text{ in } L^{p'} (0,T;L^{p'} (\Omega)) \text{ weakly};\\
\label{Emu} \partial_t v^\mu \rightharpoonup \partial_t v & \text{ in }  L^{p'} (0,T; W^{-1,p'}(\Omega)) + L^p(0,T;L^p(\O))  \text{ weakly}.
\end{align}
Note that \eqref{Bmu}, \eqref{Cmu}, \eqref{Emu} and Proposition \ref{con}
imply that
\begin{align}\label{c2a}
  v \in C([0,T];L^2(\O))\,.
\end{align}
 \noindent{\bf Limiting process} - If we use $\psi \in C^\infty_0(\Omega_T)$ as test function in  the weak formulation of $(\ref{PFep})$ in Definition $\ref{wmualpha}$, we have
 \be \label{weakff1} \int_0^T  (v^\mu, \partial_\tau \psi) -(a(\mu, v^\mu) \nabla v^\mu, \nabla \psi)   \, d\tau \; =  \int_0^T (J_\mu(v^\mu)\cdot \nabla v^\mu, \psi) \, d\tau.\ee
 Now,  as $\mu \to 0^+$, by $(\ref{Amu})$
 $$\int_0^T (v^\mu, \partial_\tau \psi) \, d \tau \to \int_0^T (v, \partial_\tau \psi) \, d \tau; $$
  by $(\ref{Dmu})$
  $$\int_0^T(a(\mu, v^\mu) \nabla v^\mu, \nabla \psi)   \, d\tau \to \int_0^T(\chi, \nabla \psi)   \, d\tau;$$
  by $(\ref{Amu})$, $(\ref{Cmu})$, and $(\ref{M1})$
  $$\int_0^T (J_\mu(v^\mu)\cdot \nabla v^\mu, \psi) \, d\tau \to \int_0^T (v\cdot \nabla v, \psi) \, d\tau.$$
  So, we get $v\in L^{p}(0,T;W^{1,p}_0(\Omega))\cap L^{\infty}(\Omega_T)$, with 
  \begin{equation}
  \partial_\tau v =  \nabla \cdot \chi - v\cdot \nabla v \in L^{p'}(0,T;W^{-1,p'}(\Omega)) + L^{p}(0,T;L^{p}(\Omega)),
  \label{vt2}
    \end{equation}
   and we need to show that $\chi = \left| \nabla v \right|^{p-2}\nabla v.$
  To this end we define
\be X_\mu := \int_0^T (a(\mu, v^\mu)\nabla v^\mu - a(\mu, \varphi)\nabla \varphi, \nabla v^\mu - \nabla \varphi) \, d \tau . \label{mumon}\ee
From the properties of the  $(p,\mu)$-Laplacian it is well known
that $X_\nu \geq 0$,  for all $ \varphi \in L^{p}(0,T;W^{1,p}_0(\Omega)).$
Note that from \eqref{weakf1} and \eqref{ues} follows
  that \linebreak${\partial_t v^\mu \in  L^{p'} (0,T;
    W^{-1,p'}(\Omega)) + L^p(0,T;L^p(\O))}$, which together with
  Proposition~\ref{existencemualpha} and  Proposition \ref{ibp} yields
  \begin{align}
\int_0^T (a(\mu, v^\mu)\nabla v^\mu, \nabla v^\mu) \, d\tau   &=
-\frac{1}{2} ||v^\mu (T)||_{ L^{2}(\Omega)}^2 + \frac{1}{2}  ||v_\circ||_{ L^{2}(\Omega)}^2 \notag\\
&\quad - \int_0^T (J_\mu (v^\mu) \cdot \nabla v^\mu, v^\mu)\, d\tau. \label{lsmu}
  \end{align}
Due to $(\ref{Amu})$, $(\ref{Cmu})$, and $(\ref{M1})$ we get
\begin{equation}\label{lscmu}
\limsup_{\mu \to 0^+} \left(\int_0^T (J_\mu (v^\mu) \cdot \nabla v^\mu, v^\mu)\, d\tau \right) = \int_0^T (v \cdot \nabla v, v)\, d\tau.
\end{equation} 
Moreover, we prove that it is possible to choose the subsequence $v^\nu$ in such a way
$v^\nu(T)\rightharpoonup v(T) \text{ in }  L^{2}(\Omega)$, and
consequently we get 
\begin{equation}\label{lsTmu}\liminf_{\mu \to 0^+} ||v^\mu(T)||_{ L^{2}(\Omega)} \geq ||v(T)||_{ L^{2}(\Omega)}.
\end{equation}
To this end we note that $v^\mu \in C(\left[ 0,T\right];
L^{2}(\Omega))$, by Proposition \ref{con}, and that $v^\mu(T)$ is
uniformly bounded with respect to $\mu$ in $ L^{2}(\Omega)$
(cf. \eqref{mp2}). Thus, there exists a subsequence, $v^\mu (T)$, such that 
$v^\mu(T)\rightharpoonup \xi_T \text{ in }  L^{2}(\Omega).$ 
Using $\psi(t,x) = \phi(t) \omega(x) \in C^\infty([0,T];
C^\infty_0(\Omega)),$ with $\phi (0)=1, \phi (T)=0$ or $\phi (0)=0,
\phi (T)=1$, as test function in \eqref{weakff1} and passing to the
limit with respect to $\mu \to 0^+$, we obtain
\begin{align*}
&\int_0^T (v, \partial_\tau \psi) \, d\tau - \int_0^T(\chi, \nabla \psi) \, d\tau \\ &= \phi(T) (\xi_T, \omega(x)) - \phi(0)(v_\circ, \omega(x)) + \int_0^T( J_\mu (v) \cdot \nabla v, \psi) \, d\tau.
\end{align*}
On the other hand, using $(\ref{vt})$ together with Proposition $\ref{ibp}$, we get 
\begin{align*}
&\int_0^T (v, \partial_\tau \psi) \, d\tau - \int_0^T(\chi, \nabla \psi) \, d\tau  \\&= \phi(T) (v(T), \omega(x)) - \phi(0)(v(0), \omega(x)) + \int_0^T( J_\mu (v) \cdot \nabla v, \psi) \, d\tau.
\end{align*}
So, we deduce that $\xi_T=v(T)$ in $ L^{2}(\Omega)$, and we get
$(\ref{lsTmu})$. Moreover, we have $v_\circ = v(0)$ in $ L^{2}(\Omega).$
Coming back to $(\ref{lsmu})$, using $(\ref{lscmu})$, $(\ref{lsTmu})$, we have
\begin{align*}
&\limsup_{\mu \to 0^+}\left( \int_0^T (a(\mu, v^\mu)\nabla v^\mu, \nabla v^\mu) \, d\tau \right) \\&\leq -\frac{1}{2}||v(T)||_{ L^{2}(\Omega)}^2 + \frac{1}{2}||v_\circ||_{ L^{2}(\Omega)}^2 - \int_0^T(v\cdot \nabla v, v)\, d\tau.
\end{align*}
Using that $a(\mu, \phi) \le |\nabla \phi|^{p-2} $ and the
  Lebesgue dominated convergence theorem we get for any $\phi \in
  L^{p}(0,T;W^{1,p}_0(\Omega))$ that 
  \begin{align*}
    a(\mu,\phi) \nabla \phi \to |\nabla \phi|^{p-2} \nabla \phi
    \textrm{ in } L^{p'}(\Omega_T) \textrm{ strongly}.
  \end{align*}
Thus, from \eqref{mumon} follows, also using \eqref{Cmu}, \eqref{Dmu}, 
\begin{align*}
 0 \leq  \limsup_{\mu \to 0^+} X_\mu  &\leq  -\frac{1}{2}||v(T)||_{ L^{2}(\Omega)}^2 + \frac{1}{2}||v_\circ||_{ L^{2}(\Omega)}^2 - \int_0^T(v\cdot \nabla v, v)\, d\tau\\ 
&\quad -\int_0^T(\chi, \nabla \varphi)\,d\tau -\int_0^T ( \left|\nabla \varphi \right|^{p-2}\nabla \varphi, \nabla v - \nabla \varphi) \, d \tau,
\end{align*}
and, since we have proved $v(0)=v_\circ$, Proposition $\ref{ibp}$
guaranties 
$$-\frac{1}{2}||v(T)||_{ L^{2}(\Omega)}^2 + \frac{1}{2}||v_\circ||_{ L^{2}(\Omega)}^2 \;  =  -\int_0^T \langle\partial_\tau v , v\rangle _{W^{1,p}_0(\Omega) \cap L^{p'}(\Omega)} \, d\tau. $$
This implies
$$0 \leq \int_0^T (\chi - \left|\nabla \varphi \right|)^{p-2}\nabla \varphi, \nabla v - \nabla \varphi) \, d\tau, \; \;\text{for all } \varphi \in L^{p}(0,T; W^{1,p}_0(\Omega)),$$
so, by Minty Trick, $\chi = a(\mu,v)\nabla v.$\\
\textbf{Maximum principle } -
In view of (\ref{mp2})  there exists for each fixed $t \in [0,T]$ a
subsequence $v^{\mu '}(t)$ of $v^{\mu }(t)$ such that $v^{\mu
  '}(t)\overset{\ast}{\rightharpoonup} \xi_t$ in
$L^{\infty}(\Omega)$. Moreover, the uniform boundedness of $v^\mu (t)$
in $C([0,T]; L^{2}(\Omega))$ (cf.~\eqref{ues}) yields that there exists
a subsequence $v^{\mu ''}(t)$ of $v^{\mu '}(t)$ such that $v^{\mu
  ''}(t)\rightharpoonup \tilde{\xi}_t$ in $ L^{2}(\Omega)$. By
uniqueness of the weak limits we have $\tilde{\xi}_t= \xi_t$. In view
of \eqref{c2a} we also have that ${v \in C([0,T]; L^{2}(\Omega))}$, and thus,
we want to show that $\xi_t = v(t)$ in $L^2(\O)$.
To this end we consider $\tilde{\psi}(\tau,x)= \phi(\tau) \omega(x)
\in C^\infty ([0,T];C^\infty_0(\Omega))$, with $\phi(t)=1$ and
$\phi(0)=0$, and define $\psi(\tau,x) := \tilde{\psi}(\tau,x)
\chi_{\left[0,t\right]} (\tau)$. Choosing $\psi=v^{\mu '}(t)$ in the
weak formulation of Definition $\ref{wmualpha} $ for $\mu'$, and passing to the limit, we get
 $$\int_0^t ( v, \partial_\tau \tilde{\psi}) \, d\tau -\int_0^t ( \left|\nabla v\right|^{p-2} \nabla v, \nabla \tilde{\psi})) \, d\tau\, = \int_0^t (v\cdot \nabla v, \tilde{\psi}) \, d\tau - (\xi_t, w).$$
 Using \eqref{vt2} and Proposition \ref{ibp}, we get
$$\int_0^t ( v, \partial_\tau \tilde{\psi}) \, d\tau -\int_0^t ( \left|\nabla v\right|^{p-2} \nabla v, \nabla \tilde{\psi})) \, d\tau\, = \int_0^t (v\cdot \nabla v, \tilde{\psi}) \, d\tau - (v(t), w),$$
which implies $\xi_t = v(t)$ in $L^2(\O)$ and
$||v(t)||_{L^{\infty}(\Omega)} \leq
||v_\circ||_{L^{\infty}(\Omega)}$, by (\ref{mp2}) and the semi-continuity of $L^{\infty}$-norm.
Using the same argument for each $t \in \left[0,T\right]$, we get ($\ref{mpf}$).
\\ \textbf{Initial datum} -
At first, we know  $$\lim_{t\to 0^+} ||v(t)-v_\circ||_{ L^{2}(\Omega)}=0.$$
  Indeed, since Proposition \ref{con}, we have $v \in C([0,T]; L^{2}(\Omega))$ and we know $v(0)=v_\circ (x) \in  L^{2}(\Omega)$. 
  Moreover, by the maximum principle, we get $(\ref{id})$ and $v\in C([0,T]; L^q(\O))$, for any $q\in [1,\infty)$.
  \setcounter{equation}{0}
\section*{\normalsize A. Appendix}
\renewcommand{\theequation}{A.\arabic{equation}}
\renewcommand{\thetho}{A.\arabic{tho}}
\renewcommand{\thedefi}{A.\arabic{defi}}
\renewcommand{\therem}{A.\arabic{rem}}
\renewcommand{\theprop}{A.\arabic{prop}}
\renewcommand{\thelemma}{A.\arabic{lemma}}
\renewcommand{\thecoro}{A.\arabic{coro}}
\setcounter{tho}{0}
We define the space
{\small $$W:= \{ u \in L^{\infty}(\Omega_T) \cap L^{p}(0,T; W^{1,p}_0(\Omega)) \,\big| \,\partial_t u \in L^{p}(0,T; L^{p}(\Omega)) + L^{p'}(0,T;W^{-1,p'}(\Omega))\} ,$$}%
endowed with the  norm
$$||u||_{W} := ||u||_{L^{\infty}(\Omega_T)\cap L^{p}(0,T; W^{1,p}_0(\Omega))} + ||\partial_t u||_{L^{p}(0,T; L^{p}(\Omega)) + L^{p'}(0,T;W^{-1,p'}(\Omega))},$$
where
$$ ||u||_{L^{\infty}(\Omega_T)\cap L^{p}(0,T; W^{1,p}_0(\Omega))} := ||u||_{L^{\infty}(\Omega_T)} + ||u||_{ L^{p}(0,T; W^{1,p}_0(\Omega))}, $$
and
\begin{equation}\label{ns}
{ \footnotesize ||\partial_t u||_{L^{p}(0,T; L^{p}(\Omega)) + L^{p'}(0,T;W^{-1,p'}(\Omega))} := \inf_{(f,g) \in \mathcal{F}(\partial_t u)} \Set{ \left|\left|f\right| \right|_{L^{p}(0,T;L^{p}(\Omega)) }+\left|\left|g\right|\right|_{L^{p'}(0,T;W^{-1,p'}(\Omega))}}\! ,}
 \end{equation}
with
$$\mathcal{F}(\partial_t u):= \Set{(f,g) | f \in L^{p}(0,T;L^{p}(\Omega)), \, g \in L^{p'}(0,T;W^{-1,p'}(\Omega)), \, f+g =\partial_t u }\!.$$
The space $W$ is embedded in the space
{\small $$\widehat{W} := \{ u \in L^{p'}(0,T;L^{p'}(\Omega)) \cap L^{p}(0,T; W^{1,p}_0(\Omega)) \,:\, \partial_t u \in L^{p}(0,T; L^{p}(\Omega)) + L^{p'}(0,T;W^{-1,p'}(\Omega))\} ,$$}%
endowed with the norm
$$||u||_{\widehat{W}} := ||u||_{L^{p'}(0,T;L^{p'}(\Omega))\cap L^{p}(0,T; W^{1,p}_0(\Omega))} + ||\partial_t u||_{L^{p}(0,T; L^{p}(\Omega)) + L^{p'}(0,T;W^{-1,p'}(\Omega))},$$
where
$$ ||u||_{L^{p'}(0,T;L^{p'}(\Omega))\cap L^{p}(0,T; W^{1,p}_0(\Omega))} := ||u||_{L^{p'}(0,T;L^{p'}(\Omega))} + ||u||_{ L^{p}(0,T; W^{1,p}_0(\Omega))}. $$
\subsubsection*{\normalsize A.1. Density result}
\begin{prop}\label{den}
For $p\in (1,2)$, $T\in\R^+$ and a bounded domain $\Omega$ there holds
$$\overline{C^\infty([0,T]; L^{p'}(\Omega) \cap W^{1,p}_0(\Omega))}_{|| \cdot ||_{\widehat{W}}} = \widehat{W}.$$
\end{prop}
\Pr If we take $u \in \widehat{W}$, we have $u \in L^{p'}(0,T;L^{p'}(\Omega)) \cap L^{p}(0,T; W^{1,p}_0(\Omega))$, so if we consider a sequence of mollifier in time $\rho_n(t)$, we get 
$$\rho_n(t) * u(t,x) \in C^\infty ([0,T]; L^{p'}(\Omega) \cap W^{1,p}_0 (\Omega))$$
and
$$\rho_n(t) * u(t,x) \to u(t,x) \text{ in } L^{p'}(0,T;L^{p'}(\Omega)), \text{ as } n \to + \infty; $$
$$\rho_n(t) * u(t,x) \to u(t,x) \text{ in } L^{p}(0,T; W^{1,p}_0(\Omega)), \text{ as } n \to + \infty. $$
Therefore, we get
\begin{align}
&||\rho_n(t)*u(t,x) -u(t,x) ||_{L^{p'}(0,T;L^{p'}(\Omega))\cap L^{p}(0,T; W^{1,p}_0(\Omega))} = \nonumber\\ 
&= ||\rho_n(t)*u(t,x) -u(t,x) ||_{L^{p'}(0,T;L^{p'}(\Omega))} \nonumber\\
& \; \; \; \;  + ||\rho_n(t)*u(t,x) -u(t,x) ||_{ L^{p}(0,T; W^{1,p}_0(\Omega))} \to 0, \text{ as } n \to +\infty.\label{convappendix}
\end{align}
Moreover, for $\partial_t u \in L^{p}(0,T; L^{p}(\Omega)) +
L^{p'}(0,T;W^{-1,p'}(\Omega)),$ assume that \linebreak ${u^1 \in L^{p}(0,T; L^{p}(\Omega))}$ and $u^2 \in L^{p'}(0,T;W^{-1,p'}(\Omega))$ are such that $$\partial_t u = u^1+u^2.$$
We observe 
$$(\rho_n(t) * \partial_t u(t,x)) = \rho_n(t) * (u^1(t,x)) + \rho_n(t) * (u^2(t,x)). $$
We have
\begin{align*}
\rho_n(t) * (u^1(t,x)) \to u^1(t,x) &\text{ in } L^{p}(0,T;L^{p}(\Omega)), \text{ as } n \to + \infty;\\
\rho_n(t) * (u^2(t,x)) \to u^2(t,x) &\text{ in } L^{p'}(0,T;W^{-1,p'}(\Omega)), \text{ as } n \to + \infty,
\end{align*}
and thus, we get
\begin{align}
&||\rho_n(t)* \partial_t u(t,x) - \partial_t u(t,x)||_{L^{p}(0,T; L^{p}(\Omega)) + L^{p'}(0,T;W^{-1,p'}(\Omega))} \nonumber\\
&\leq||\rho_n(t)*u^1(t,x) - u^1(t,x)||_{L^{p}(0,T; L^{p}(\Omega)) + L^{p'}(0,T;W^{-1,p'}(\Omega))} \nonumber \\
& \; \; \; \;+ ||\rho_n(t)*u^2(t,x) - u^2(t,x)||_{L^{p}(0,T; L^{p}(\Omega)) + L^{p'}(0,T;W^{-1,p'}(\Omega))} \nonumber \\
&\leq||\rho_n(t)*u^1(t,x) - u^1(t,x)||_{L^{p}(0,T; L^{p}(\Omega))}\nonumber \\
&\; \; \; \; + ||\rho_n(t)*u^2(t,x) - u^2(t,x)||_{ L^{p'}(0,T;W^{-1,p'}(\Omega))} \to 0, \text{ as } n \to +\infty, \label{conv2appendix}
\end{align}
where the last inequality is due to (\ref{ns}).
Ultimately, \eqref{convappendix} and \eqref{conv2appendix} yield
$$||\rho_n(t)*u(t,x) - u(t,x)||_{\widehat{W}} \to 0, \text{ as } n \to + \infty.$$
\chiu
\subsubsection*{\normalsize A.2 Continuity}
\begin{prop}\label{con}
For $p\in (1,2)$, $T\in\R^+$ and a bounded domain $\Omega$, there holds
$$W \subset C([0,T]; L^{2}(\Omega)).$$
\end{prop}
\Pr By Theorem 1.1 in \cite{Por}, we get
$$W \subset C([0,T];L^1(\Omega)).$$
Moreover, since $u \in L^{\infty} (\Omega_T),$ we have 
$u \in L^{\infty}(0,T;L^q(\Omega))$, $q \in
[1, + \infty)$. 
So, $u\in W$ implies $u \in C([0,T];L^1(\Omega)) \cap L^{\infty}(0,T; L^q(\Omega))$. Thus, for fixed $q>2$ in \cite[Theorem 2.1]{Str} or \cite[Lemma 3.1.4]{Tem}, gives
\begin{equation} \label{Temam}
u\in C_w([0,T];L^q(\Omega)) \; \text{ and } \; ||u(t)||_{L^q (\Omega)} \leq M, \text{ for all } t\in [0,T].
\end{equation}
Now, we get $u \in C([0,T]; L^{2}(\Omega))$, since for each $t_1, t_2 \in [0,T]$
\begin{align*}
||u(t_1)-u(t_2)||_{ L^{2}(\Omega)} &\leq ||u(t_1)-u(t_2)||^a_{L^1(\Omega)} ||u(t_1)-u(t_2)||^{1-a}_{L^q(\Omega)}\\
& \leq (2 M)^{1-a} ||u(t_1)-u(t_2)||^{a}_{L^1(\Omega)},
\end{align*}
with $\frac{1}{2}= \frac{a}{1} + \frac{1-a}{q}$.
\chiu
\subsubsection*{\normalsize A.3 Integration by parts formula}
\begin{prop}\label{ibp}
For $p\in (1,2)$, $T\in\R^+$ and a bounded domain $\Omega$ there holds
for all $u, v \in W $ that
\begin{align*}
&\int_0^t\langle u,\partial_\tau v\rangle _{L^{p'}(\Omega)\cap W^{1,p}_0(\Omega)} \, d\tau \\
  &=
 (u(t),v(t)) - (u(0),v(0)) 
 - \int_0^t\langle v, \partial_\tau u\rangle _{L^{p'}(\Omega)\cap W^{1,p}_0(\Omega) } \, d\tau,
\end{align*}
for all $t \in \left[0,T\right].$
\end{prop}
\Pr
Let us consider $u,v \in W$. Then Proposition \ref{con} implies that  $u,v \in C([0,T]; L^{2}(\Omega))$, and that $(u(t),v(t))\in \R$ is well defined  for all $ t \in \left[0,T\right].$
Due to Proposition \ref{den}, we know that there exist  $$\lbrace u_n \rbrace_{n \in \N}, \lbrace v_n \rbrace_{n \in \N} \in C^\infty ([0,T];W^{1,p}_0(\Omega)\cap L^{p'}(\Omega))$$
such that
$$u_n \to u \text{ in } \widehat{W};$$
$$v_n \to v \text{ in } \widehat{W}.$$
In particular,
\begin{equation}\label{fu}
u_n \to u \text{ in } L^{p}(0,T;W^{1,p}_0(\Omega)) \cap L^{p'}(0,T;L^{p'}(\Omega)).
\end{equation}
Moreover, if we  consider $u^1\in L^{p}(0,T;L^{p}(\Omega))$ and $u^2 \in L^{p'}(0,T;W^{-1,p'}(\Omega))$ such that $\partial_t u = u^1+u^2$, we know, for the sequence built as in Proposition \ref{den}, that $\partial_t (u_n)= (\partial_t u)_n =(u^1)_n + (u^2)_n  $ and
\begin{align}\label{fu1}(u^1)_n \to u^1 &\text{ in } L^{p}(0,T;L^{p}(\Omega)), \\
\label{fu2} (u^2)_n \to u^2 &\text{ in } L^{p'}(0,T;W^{-1,p'}(\Omega)).
\end{align}
The same is true for $v_n$ and $v$:
\begin{equation}\label{fv}
v_n \to v \text{ in } L^{p}(0,T;W^{1,p}_0(\Omega)) \cap L^{p'}(0,T;L^{p'}(\Omega)).
\end{equation}
Moreover, if we  consider $v^1\in L^{p}(0,T;L^{p}(\Omega))$ and $v^2 \in L^{p'}(0,T;W^{-1,p'}(\Omega))$ such that $\partial_t v = v^1+v^2$, we know, for the sequence built as in Proposition \ref{den}, that $\partial_t (v_n) = (\partial_t v)_n = (v^1)_n + (v^2)_n  $
\begin{align}\label{fv1}(v^1)_n \to v^1 &\text{ in } L^{p}(0,T;L^{p}(\Omega)), \\
\label{fv2} (v^2)_n \to v^2 &\text{ in } L^{p'}(0,T;W^{-1,p'}(\Omega)).
\end{align}
So, since (\ref{fu}), (\ref{fv}),
and  Hölder inequaity $(u_n(\tau), v_n(\tau)) \to (u(\tau),v(\tau))$ in $L^1 (0,T;\R)$, so 
\begin{equation}\label{dis}
\frac{d}{d\tau} (u_n(\tau), v_n(\tau)) \to \frac{d}{d\tau} (u(\tau),v(\tau)) \text{ in } D' (0,T;\R),
\end{equation}
where $D' (0,T;\R)$ is the space of distribution.
For shortage of notation we  define
$$\langle\cdot,\cdot\rangle _{L^{p'}(\Omega)\cap W^{1,p}_0(\Omega)} =: \langle\cdot ,\cdot\rangle .$$
Since $(u_n(\tau),v_n(\tau)) \in C^\infty(0,T;\R)$ and $(\cdot,\cdot)$ is bilinear and continuous, we get
$$
 \frac{d}{d\tau} (u_n(\tau),v_n(\tau)) =
 \langle v_n(\tau),\partial_\tau (u_n)(\tau)\rangle + \langle u_n(\tau),\partial_\tau (v_n)(\tau)\rangle ,
$$
and we  show that
\begin{align}\label{q}
\langle v_n(\tau),\partial_\tau (u_n)(\tau)\rangle  \to \langle v(\tau), \partial_\tau u(\tau)\rangle & \text{ in } L^1(0,T;\R)\\
\label{w} \langle u_n(\tau), \partial_\tau (v_n)(\tau)\rangle  \to\langle u(\tau),\partial_\tau v(\tau)\rangle &  \text{ in } L^1(0,T;\R).
\end{align} 
In detail for (\ref{q}), but for (\ref{w}) it is the same, we have
\begin{align*}
\langle v_n(\tau), \partial_\tau(u_n)(\tau)\rangle  & = \langle v_n(\tau), (u^1 +u^2 )_n(\tau)  \rangle  \\
& =   \langle v_n(\tau), (u^1)_n(\tau)\rangle  + \langle v_n(\tau), (u^2)_n(\tau)\rangle ,
\end{align*}
for the first term of the sum we use (\ref{fu1}) and (\ref{fv}), for the second one (\ref{fu2}) and (\ref{fv}).
So, using (\ref{dis}), (\ref{q}), (\ref{w}), we get
$$\frac{d}{d\tau }(u(\tau), v(\tau))  = \langle u(\tau), \partial_\tau
v(\tau)\rangle  + \langle v(\tau), \partial_\tau u(\tau)\rangle \in L^1(0,T;\R),$$
therefore 
$$(u(\tau), v(\tau)) \in W^{1,1}(0;T;\R)$$
and we  write it as the integral of its derivative:
\begin{align*}
(u(t), v(t)) - (u(0), v(0)) &= \int_0^t \frac{d}{d\tau} (u(\tau), v(\tau)) \, d\tau \\
&= \int_0^t \left( \langle u(\tau), \partial_\tau v(\tau)\rangle  + \langle v(\tau), \partial_\tau u(\tau)\rangle  \right) \, d\tau.\end{align*}
\chiu\\
 {\bf Acknowledgment} -
The research activity of A.D. is performed under the auspices of the group GNFM of INdAM.
\vskip0.1cm\noindent
 {\bf Declarations}
\vskip0.1cm\noindent
 {\bf Funding} -  A.D. acknowledges the Institute of Applied Mathematics, Albert-Ludwigs-University, Freiburg, for the support provided during the period of work on this paper.
\vskip0.1cm\noindent
 {\bf Conflict of interest} - The authors have no conflicts of interest to declare that are relevant to the content of this article.

\end{document}